\documentclass[twoside,a4paper]{article}

\usepackage{amsmath,amsthm,amssymb,latexsym}
\usepackage[v2,tips]{xy}

\theoremstyle{plain}%
\newtheorem{proposition}{Proposition}[section]%
\newtheorem{theorem}[proposition]{Theorem}%
\newtheorem{lemma}[proposition]{Lemma}%
\theoremstyle{definition}%
\newtheorem{definition}[proposition]{Definition}%

\newcommand{\vnten}{\bar\otimes}

\newcommand{\proten}{{\widehat{\otimes}}}
\newcommand{\mc}[1]{\mathcal{#1}}
\newcommand{\ip}[2]{{\langle {#1} , {#2} \rangle}}
\newcommand{\inten}{{\check{\otimes}}}

\begin{document}

\large
\title{\textsc{$p$-Operator Spaces and Fig\`a-Talamanca-Herz Algebras}}
\author{Matthew Daws}
\maketitle

\begin{abstract}
We study a generalisation of operator spaces modelled on $L_p$ spaces,
instead of Hilbert spaces, using the notion of $p$-complete boundedness,
as studied by Pisier and Le Merdy.  We show that the Fig\`a-Talamanca-Herz
Algebras $A_p(G)$ becomes quantised Banach algebras in this framework,
and that the cohomological notion of amenability of these algebras
corresponds to amenability of the locally compact group $G$.  We thus
argue that we have presented a generalised of the use of operator spaces
in studying the Fourier algebra $A(G)$, in the spirit of Ruan.  Finally,
we show that various notions of multipliers of $A_p(G)$ (including Herz's
generalisation of the Fourier-Stieltjes algebra) naturally fit into
this framework.


\end{abstract}

\section{Introduction}

The Fourier algebra, $A(G)$, of a locally compact group $G$ is the
collection of coefficient functionals $f:G\rightarrow\mathbb C$
of the form
\[ f(g) = [\lambda(g)(x),y] \qquad (g\in G), \]
where $x,y\in L_2(G)$ and $\lambda$ is the left-regular representation
of $G$ on $L_2(G)$.  Eymard defined and studied this commutative
Banach algebra in \cite{Eymard}.  For an abelian group $G$, the
Fourier transform shows that $A(G)$ is nothing by $L_1(\hat G)$,
where $\hat G$ is the dual group of $G$.  As such, $A(G)$ is
amenable as a Banach algebra, and for another abelian group $H$,
we have that $A(G)\proten A(H) = A(G\times H)$.  However, as first
noted by Johnson in \cite{J1}, there exist compact groups $G$
for which $A(G)$ is not amenable.  Thus the Banach algebra $A(G)$
does not seem to capture some properties of the group $G$.

In \cite{Ruan1}, Ruan showed that when $A(G)$ is considered
as an operator space (and hence as a quantised Banach algebra),
we have that $A(G)$ is amenable if and only if $G$ is amenable,
and that $A(G)\proten A(H) = A(G\times H)$ for all locally compact
groups $G$ and $H$ (here we use the operator space projective tensor
product).  These results provide some compelling evidence that
$A(G)$ is best viewed as an operator space, and not simply as a
Banach algebra.

In \cite{Fig}, Fig\`a-Talamanca introduced a natural generalisation
of the Fourier algebra, for abelian and compact groups, by replacing $L_2(G)$
by $L_p(G)$.  In \cite{Herz1}, Herz extended the definition to arbitrary groups,
leading to the commutative Banach algebra $A_p(G)$, now called
the Fig\`a-Talamanca-Herz algebras.  In many ways these algebras
behave like $A(G)$; for example, Leptin's theorem (see
\cite[Theorem~6]{Herz2} or \cite[Section~10]{Pier}) states that
$G$ is an amenable group if and only if $A_p(G)$ has a
bounded approximate identity.

There have been a number of attempts to give $A_p(G)$ an operator
space structure.  In \cite{Runde1}, Runde used some of Pisier's work
on interpolation spaces to define an operator space version of $A_p(G)$,
denoted $OA_p(G)$.  Unfortunately, while $OA_2(G)=A(G)$ as Banach spaces,
the operator space structure can differ; furthermore, $OA_p(G)$ can
fail to be equal to $A_p(G)$, even as a Banach space, for $p\not=2$.
In \cite{LNR}, the authors use Lambert's ideas of row and column
operator spaces to define an operator space structure on $A_p(G)$
which turns $A_p(G)$ into a bounded (but not contractive) quantised
Banach algebra, and in such a way that $A_2(G) = A(G)$ completely
isometrically.  Furthermore, $A_p(G)$ is amenable in this framework
if and only if $G$ is an amenable group.

In this paper, we shall use ideas of Pisier and Le Merdy to define
the notion of a \emph{$p$-operator space} (for $1<p<\infty$, with a
$2$-operator space being simply an operator space).  We show that
the algebras $A_p(G)$ then carry a natural $p$-operator space
structure.  We investigate the amenability of $A_p(G)$ in this
framework, and also study the $p$-completely bounded multipliers
of $A_p(G)$.

\section{Banach spaces}

In this section we shall gather together some basic results on
Banach spaces.  Let $E$ be a Banach space, and denote by $E'$ the
dual space of $E$.  For $x\in E$ and $\mu\in E'$, we write
$\ip{\mu}{x}$ for $\mu(x)$ (we use angle brackets for bilinear
products, and occasionally use square brackets for sesquilinear
products).  There is a canonical isometry $\kappa_E:E\rightarrow E''$
defined by $\ip{\kappa_E(x)}{\mu} = \ip{\mu}{x}$.  When $\kappa_E$ is
an isomorphism, we say that $E$ is \emph{reflexive}

Let $E$ and $F$ be Banach spaces, and consider
the algebraic tensor product $E\otimes F$.  We define the
\emph{projective tensor norm} $\|\cdot\|_\pi$ on $E\otimes F$ by
\[ \|\tau\|_\pi = \inf\Big\{ \sum_k \|x_k\| \|y_k\| :
\tau=\sum_k x_k \otimes y_k \Big\} \qquad (\tau\in E\otimes F). \]
The completion of $E\otimes F$ with respect to $\|\cdot\|_\pi$ is
denoted by $E\proten F$.  It is a simple exercise to show that
$(E\proten F)' = \mc B(E,F') \cong \mc B(F,E')$ by the identification
\[ \ip{T}{x\otimes y} = \ip{T(x)}{y} \qquad
(T\in\mc B(E,F'), x\in E, y\in F). \]
Here we write $\mc B(E,F)$ for the Banach space of bounded linear
operators from $E$ to $F$.  We write $\mc B(E)$ for $\mc B(E,E)$.

Alternatively, we may embed $E\otimes F$ into $\mc B(E',F)$, which
leads to the definition of the \emph{injective tensor norm}
$\|\cdot\|_\epsilon$, and the \emph{injective tensor product}
$E\inten F$.  Then $E'\otimes F$ can be identified with the
\emph{finite rank operators} from $E$ to $F$, denoted by $\mc F(E,F)$.
The closure of $\mc F(E,F)$ in $\mc B(E,F)$ is the \emph{approximable
operators} from $E$ to $F$, denoted by $\mc A(E,F)$.  Thus $E'\inten F
= \mc A(E,F)$.

There is an obvious norm-decreasing map $J:E'\proten E\rightarrow
E'\inten E = \mc A(E)$, whose image is the \emph{nuclear operators},
$\mc N(E)$.  We give $\mc N(E)$ the quotient norm coming from
$\mc N(E) \cong E' \proten E / \ker J$.  When $J$ is injective,
we say that $E$ has the \emph{approximation property}.
See \cite{Ryan} or \cite{DF} for further details on these ideas.

\section{Amenable Banach algebras}\label{amen_ba}

We shall eventually apply our results to the study of when certain
Banach algebras are \emph{amenable} (in various senses).  However,
we shall also need some ideas from this area as we go along, so we
introduce the needed ideas now.

Let $\mc A$ be a Banach algebra, and let $E$ be an $\mc A$-bimodule.
A linear map $d:\mc A\rightarrow E$ is a \emph{derivation} if
$d(ab) = a\cdot d(b) + d(a)\cdot b$ for $a,b\in\mc A$.  We shall assume
that all our derivations are bounded.  For $x\in E$, define
$d_x:\mc A\rightarrow E$ by $d_x(a)=a\cdot x-x\cdot a$, for $a\in \mc A$.
Then $d_x$ is a derivation, called an \emph{inner derivation}.
A Banach algebra $\mc A$ is \emph{amenable} when every derivation
from $\mc A$ to a dual bimodule is inner.  See the book \cite{RundeBook}
for details about amenable Banach algebras, for example.

Johnson showed in \cite{Johnson} that for a locally compact group $G$,
one has that $G$ is amenable if and only if the group algebra $L_1(G)$
is amenable.  Recall that a group $G$ is \emph{amenable} when there
is a left-invariant mean for $L_\infty(G)$.  See \cite{Paterson} or
\cite{Pier} for details about amenable groups.
Johnson also provided a useful characterisation of when
an algebra is amenable.

\begin{definition}
Let $\mc A$ be a Banach algebra.  A bounded net $(d_\alpha)$ in
$\mc A\proten\mc A$ is an \emph{approximate diagonal} if
\[ \lim_\alpha \|a\cdot d_\alpha-d_\alpha\cdot a\|=0,
\quad \lim_\alpha \|a\Delta_{\mc A}(d_\alpha)-a\|=0 \qquad (a\in\mc A). \]
Here $\Delta_{\mc A}:\mc A\proten\mc A\rightarrow\mc A$ is the
linearisation of the product, defined by $\Delta_{\mc A}(a\otimes b)=ab$.
\end{definition}

\begin{theorem}\label{amen_char}
Let $\mc A$ be a Banach algebra.  Then $\mc A$ is amenable if
and only if $\mc A$ has an approximate diagonal.  When $G$ is an
amenable group, we may choose an approximate diagonal for $L_1(G)$
which is bounded by $1$.
\end{theorem}
\begin{proof}
See \cite[Chapter~2]{RundeBook} for example.
\end{proof}

Let $\mc A$ be a Banach algebra, and let $E$ be a left $\mc A$-module.
Let $\mc A_E^c = \{ T\in\mc B(E) : T(a\cdot x)=a\cdot T(x)
\ (a\in\mc A, x\in E) \}$, the commutant of $\mc A$ in $E$.  Then
a projection $\mc Q:\mc B(E)\rightarrow\mc A_E^c$ is a \emph{quasi-expectation}
when $\mc Q(TSR) = T \mc Q(S) R$ for $T,R\in\mc A_E^c$ and $S\in\mc B(E)$.

\begin{proposition}\label{quasi_exp}
Let $\mc A$ be an amenable Banach algebra, and let $E$ be a 
reflexive left $\mc A$-module.  Then there is a quasi-expectation
$\mc Q:\mc B(E)\rightarrow \mc A_E^c$.
\end{proposition}
\begin{proof}
We sketch a proof (see \cite[Theorem~4.4.11]{RundeBook} for example).
Let $(d_\alpha)$ be an approximate diagonal for $\mc A$, and
let $d_\alpha = \sum_{n=1}^\infty a^{(\alpha)}_n \otimes b^{(\alpha)}_n$
for each $\alpha$.  As $E$ is reflexive, by moving to a subnet if
necessary, we may define
\[ \ip{\mu}{\mc Q(T)(x)} = \lim_\alpha \sum_{n=1}^\infty
\ip{\mu}{a^{(\alpha)}_n\cdot T(b^{(\alpha)}_n\cdot x)}
\qquad (x\in E,\mu\in E',T\in\mc B(E)). \]
Then $\mc Q$ is a linear operator, and $\|\mc Q\|\leq\limsup_\alpha
\|d_\alpha\|$.  Clearly, if $T\in\mc A^c_E$, then $\mc Q(T)=T$.
Conversely, as $d_\alpha$ is an approximate diagonal,
for $x\in E,\mu\in E',a\in\mc A$ and $T\in\mc B(E)$,
\[ \ip{\mu}{\mc Q(T)(a\cdot x) - a\cdot\mc Q(T)(x)}
= \lim_\alpha \sum_{n=1}^\infty
\ip{\mu}{a^{(\alpha)}_n\cdot T(b^{(\alpha)}_na\cdot x)
- aa^{(\alpha)}_n\cdot T(b^{(\alpha)}_n\cdot x)} = 0, \]
so that $\mc Q(T)\in \mc A_E^c$.  Thus $\mc Q$ is a projection
onto $\mc A_E^c$.  Similarly, for $T,R\in\mc A_E^c$ and $S\in\mc B(E)$,
it is easy to check that $\mc Q(TSR) = T \mc Q(S) R$.
\end{proof}

It is shown in \cite{Daws} that, in a certain sense, the converse
to the above is true.  When $\mc A$ is a von Neumann algebra,
we follow \cite{RundeBook} and define $\mc A$ to be \emph{Connes-amenable}
using the same definition as for amenability, but insisting that everything
is suitably weak$^*$-continuous (this is commonly just referred to as
the suitable definition of ``amenable'' for von Neumann algebras).
Then $\mc A$ is Connes-amenable if and only if there is an \emph{expectation}
(that is, a norm-one projection) from $\mc B(H)$ to $\mc A$, where
$H$ is any Hilbert space such that $\mc A \subseteq \mc B(H)$ is a
concrete realisation of the von Neumann algebra $\mc A$.  It is well-known
(see \cite[Chapter~III, Theorem~3.4]{tak1}) that an expectation is always a quasi-expectation.

\section{$p$-operator spaces}

Let $SQ_p$ be the collection of subspaces of quotients of $L_p$
spaces, where we identify spaces which are isometrically isomorphic.
Let $\mu$ be a measure, and $E$ a Banach space.  We define a
norm on the algebraic tensor product $L_p(\mu)\otimes E$ by embedding
$L_p(\mu)\otimes E$ into $L_p(\mu,E)$ in the obvious way.  Let the
completion be denoted by $L_p(\mu) \otimes_p E$.  It is easy to see
that $L_p(\mu)\otimes E$ is dense in $L_p(\mu,E)$, so that
$L_p(\mu) \otimes_p E = L_p(\mu,E)$ isometrically.  An important
property of $SQ_p$ spaces is the following.  For $E,F\in SQ_p$, we have
that for $T\in\mc B(L_p(\mu))$ and $S\in\mc B(E,F)$, the operator
$T\otimes S$ is bounded as an operator from $L_p(\mu)\otimes_p E$
to $L_p(\mu)\otimes_p F$, with norm $\|T\| \|S\|$.  See
\cite[Section~7]{DF} or the survey paper \cite{DF1} for further information.

For $n\in\mathbb N$, let $\ell_p^n$ be $\mathbb C^n$ with the
$\ell_p$-norm.  Similarly, $\ell_p(I)$ is the usual $\ell_p$ space
over an index set $I$; we set $\ell_p$ to be $\ell_p(\mathbb N)$.
Throughout, we shall let $p'$ be the conjugate index to $p$,
so that $p^{-1} + p'^{-1} = 1$.

An abstract characterisation of $SQ_p$ spaces is the following, which
goes back to Kwapien (see \cite[Theorem~3.2]{Merdy} for example).
For a square matrix $a=(a_{ij})\in\mathbb M_n$, we let $a$ induce an
operator on $\ell_p^n$, which leads to the norm
\[ \|a\|_{\mc B(\ell_p^n)}
= \sup\Big\{ \Big( \sum_{i=1}^n \Big| \sum_{j=1}^n a_{ij} x_j
\Big|^p \Big)^{1/p} : (x_j)_{j=1}^n \subseteq\mathbb C,
\sum_{j=1}^n |x_j|^p \leq 1 \Big\}. \]
We have that $E\in SQ_p$ if and only if, for each $n$ and
each $a=(a_{ij})\in\mathbb M_n$, we have that
\[ \sup\Big\{ \Big(
\sum_{i=1}^n \Big\| \sum_{j=1}^n a_{ij} x_j
\Big\|^p \Big)^{1/p} : (x_j)_{j=1}^n \subseteq E,
\sum_{j=1}^n \|x_j\|^p \leq 1 \Big\} \leq \|a\|_{\mc B(\ell_p^n)}. \]

\subsection{$p$-operator spaces}

We now introduce some ideas studied in \cite{Pisier}, and especially
\cite{Merdy}, although we introduce some new notation.
A \emph{concrete $p$-operator space} is a closed subspace of
$\mc B(E)$, for some $E\in SQ_p$.  Notice that we could equally
define this by using $\mc B(E,F)$ instead, for $E,F\in SQ_p$.
This follows, as we can identify $\mc B(E,F)$ with a closed subspace
of $\mc B(E\oplus_p F)$, where $E\oplus_p F$ is the direct sum of
$E$ and $F$ together with the norm $\|e \oplus f\| = (\|e\|^p+\|f\|^p)^{1/p}$
for $e\in E$ and $f\in F$.

For a concrete $p$-operator space $X\subseteq\mc B(E)$, for each
$n>0$, we define a norm $\|\cdot\|_n$ on $\mathbb M_n(X) = \mathbb M_n
\otimes X$ by identifying $\mathbb M_n(X)$ as a subspace of
$\mc B(\ell_p^n \otimes_p E)$.  It is easy to see that the norms
$\|\cdot\|_n$ satisfy:
\newcommand{\Axone}{$\mc{D}_\infty$}
\newcommand{\Axtwo}{$\mc{M}_p$}
\begin{itemize}
\item[\Axone] for $u\in\mathbb M_n(X)$ and $v\in\mathbb M_m(X)$,
   we have that $\|u\oplus v\|_{n+m} = \max(\|u\|_n, \|v\|_m)$.  Here
   $u\oplus v\in\mathbb M_{n+m}(X)$ has block representation
   $\begin{pmatrix} u & 0 \\ 0 & v \end{pmatrix}$.
\item[\Axtwo] for $u\in\mathbb M_m(X)$, $\alpha\in\mathbb M_{n,m}$
   and $\beta\in\mathbb M_{m,n}$, we have that $\|\alpha u \beta\|_n
   \leq \|\alpha\| \|u\|_m \|\beta\|$.  Here $\alpha u \beta$ is the obvious
   matrix product, and we define $\|\alpha\|$ to be the norm of $\alpha$
   as a member of $\mc B(\ell_p^m,\ell_p^n)$, and similarly for $\beta$.
\end{itemize}

An \emph{abstract $p$-operator space} is a Banach space $X$ together with
a family of norms $\|\cdot\|_n$ defined by $\mathbb M_n(X)$ satisfying the
above two axioms.  When $p=2$, the above axioms are just Ruan's axioms,
and so $2$-operator spaces are just operator spaces.
Here, and throughout, we refer to \cite{ER} for details on operator spaces.
Then \cite[Theorem~4.1]{Merdy} shows that an abstract $p$-operator space $X$
can be isometrically embedded in $\mc B(E)$ for some $E\in SQ_p$, and in
such a way that the canonical norms on $\mathbb M_n(X)$ arising from
this embedding agree with the given norms.  Henceforth, we shall just
talk of $p$-operator spaces.  We shall tend to abuse notation, and
write $\|\cdot\|$ instead of $\|\cdot\|_n$, where there can be no
confusion.

The natural morphisms between $p$-operator spaces are the 
\emph{$p$-completely bounded} maps, as first studied in \cite{Pisier}.
A linear map $u:X\rightarrow Y$ between $p$-operator spaces induces
a map $(u)_n:\mathbb M_n(X)\rightarrow\mathbb M_n(Y)$ in an obvious
way.  We say that $u$ is \emph{$p$-completely bounded} if
$\|u\|_{pcb} := \sup_n \|(u)_n\| <\infty$.  Similarly, we have the notions
of $p$-completely contractive and $p$-completely isometric.
We write $\mc{CB}_p(X,Y)$ for the Banach space of all $p$-completely
bounded maps from $X$ to $Y$.

Pisier proved a factorisation scheme for $p$-completely bounded maps.
Let $E\in\mc SQ_p$, let $J$ be some index set, and let $\phi_j$ be
a measure, for each $j\in J$.  Let $\mc U$ be an ultrafilter on $J$, so
that we may form the ultraproduct $\hat E = (L_p(\phi_j,E))_{\mc U}$.
Notice that $\hat E\in SQ_p$ (see \cite{Hein} for details about ultraproducts
of Banach spaces).  For each $j\in J$, $\mc B(E)$ acts naturally on
$L_p(\phi_j,E)$, and so we get a canonical homomorphism $\pi:
\mc B(E)\rightarrow\mc B(\hat E)$.  Now suppose that $X\subseteq\mc B(E)$
is a $p$-operator space.  Let $N\subseteq M\subseteq\hat E$ and $\hat N
\subseteq\hat M\subseteq\hat E$ be closed subspaces such that, for each
$x\in X$, $\pi(x)$ maps $N$ into $\hat N$ and $M$ into $\hat M$.  Hence,
for each $x\in X$, $\pi(x)$ naturally induces a map, denoted $\hat\pi(x)$,
from $G = M/N$ to $\hat G=\hat M/\hat N$.  Notice that $G,\hat G\in SQ_p$.
We call the map $\hat\pi$ a \emph{$p$-representation} from $X$ to $\mc B(G,\hat G)$.

\begin{theorem}\label{pisier_pcb}
Let $E,F\in SQ_p$, let $X\subseteq\mc B(E)$ be a $p$-operator space,
and let $u:X\rightarrow\mc B(F)$ be a linear map.  Then $u$ is $p$-completely
bounded with $\|u\|_{pcb}\leq C$ if and only if there exists a $p$-representation
$\hat\pi:X\rightarrow\mc B(G,\hat G)$ and operators $U:F\rightarrow G$ and
$V:\hat G \rightarrow F$ such that
\[ u(x) = V \hat\pi(x) U \qquad (x\in X). \]
\end{theorem}
\begin{proof}
This is \cite[Theorem~2.1]{Pisier}, although we have followed the
presentation of \cite{Merdy}.
\end{proof}

As noted by Pisier after the statement of \cite[Theorem~2.1]{Pisier},
if $X\subseteq\mc B(E)$ is a unital closed subalgebra, we may
suppose that $M=\hat M$ and $N=\hat N$, so that $G=\hat G$.

As for operator spaces, we define a norm on $\mathbb M_n(\mc{CB}_p(X,Y))$
by identifying this space with $\mc{CB}_p(X,\mathbb M_n(Y))$.  It is
then an easy check to see that these norms satisfy the above axioms,
and so Le Merdy's theorem tells us that $\mc{CB}_p(X,Y)$ is itself
a $p$-operator space.

For the next result, we give $\mathbb C$ the obvious $p$-operator
space structure: that is, $\mathbb M_n(\mathbb C) = \mc B(\ell_p^n)$.

\begin{lemma}
Let $X$ be a $p$-operator space, and let $\mu\in X'$, the Banach
dual space of $X$.  Then $\mu$ is $p$-completely bounded as a map to
$\mathbb C$, and $\|\mu\|_{pcb} = \|\mu\|$.
\end{lemma}
\begin{proof}
We cannot simply follow the usual operator-space proof.  In the
$p=2$ case, we have Smith's Lemma available, which tells us that
for a map $u:X\rightarrow\mathbb M_n$, we have that $\|u\|_{cb}
= \|(u)_n\|$.  An examination of the proof of \cite[Lemma~2.2.1]{ER}
shows that we cannot hope for an extension to the general $p$ case.

We wish to show that $(\mu)_n:\mathbb M_n(X)\rightarrow\mc B(\ell_p^n)$
is bounded, with norm $\|\mu\|$.  Let $x = (x_{ij})_{i,j=1}^n \in
\mathbb M_n(X)$, so that $(\mu)_n(x) = (\ip{\mu}{x_{ij}})$.  Let
$\alpha = (\alpha_i)_{i=1}^n\in\ell_p^n$ and $\beta = (\beta_j)_{j=1}^n
\in \ell_{p'}^n$.  Then
\[ \ip{\beta}{(\mu)_n(x)(\alpha)} = \sum_{i,j=1}^n
\beta_i \ip{\mu}{x_{ij}} \alpha_j = \ip{\mu}{
\sum_{i,j=1}^n \beta_i x_{ij} \alpha_j}. \]
We may regard $\alpha$ as a member of $\mathbb M_{n,1}$, from which it
follows that $\|\alpha\|_{\mc B(\ell_p^1,\ell_p^n)} = \|\alpha\|_p$,
and similarly $\beta\in\mathbb M_{1,n}$ with $\|\beta\|_{\mc B(\ell_p^n,
\ell_p^1)} = \|\beta\|_{p'}$.  So from axiom \Axtwo\
it follows that $\|\beta x \alpha\|_1 \leq \|\beta\|_{p'} \|x\|_n
\|\alpha\|_p$, and so
\[ |\ip{\beta}{(\mu)_n(x)(\alpha)}| \leq \|\mu\|
\|\beta\|_{p'} \|x\|_n \|\alpha\|_p. \]
This implies that $\|(\mu)_n(x)\| \leq \|\mu\| \|x\|_n$, which in
turn implies that $\|(\mu)_n\| \leq \|\mu\|$, as required.
\end{proof}

As this proof indicates, we shall have significant problems extending
many results from operator spaces to $p$-operator spaces.  Indeed,
the evidence below suggests that the current definitions might be wrong,
in that we are unable to prove simple properties which one would naturally
want to hold.

We may hence identify $X'$ with $\mc{CB}_p(X,\mathbb C)$, and from
this it follows that $X'$ is also a $p$-operator space.
We may use Le Merdy's Theorem to show that $X'$ admits a representation
$X' \subseteq \mc B(E)$ for some $E\in SQ_p$.  In fact, in this special
case, we have a more concrete embedding.

\begin{theorem}\label{dual_embed}
Let $X$ be a $p$-operator space.  There exists a $p$-complete isometry
$\Phi:X'\rightarrow\mc B(\ell_p(I))$ for some index set $I$.
\end{theorem}
\begin{proof}
We follow \cite[Proposition~3.2.4]{ER}.
For each $n\in\mathbb N$, let $s_n$ be the unit sphere of $\mathbb M_n(X)$,
and let $s = \bigcup_n s_n$.  For $x\in s$, let $n(x)\in\mathbb N$ be such
that $x\in s_{n(x)}$.  Then let $E$ be the $\ell_p$-direct sum of the
spaces $\{ \ell_p^{n(x)} : x\in s\}$, so that $E$ is isometric to $\ell_p(I)$
for some index set $I$.  For $\mu\in X'$ and $x\in s$, we have that $x(\mu)
\in \mathbb M_{n(x)} = \mc B(\ell_p^{n(x)})$, with $\|x(\mu)\| \leq
\|x\| \|\mu\| = \|\mu\|$.  For $a=(a_x)_{x\in s} \in E$
and $\mu\in X'$, we may hence define
\[ \Phi(\mu)(a) = \big( x(\mu)(a_x) \big)_{x\in s}, \]
and we see that $\Phi$ is norm-decreasing.  Indeed, clearly $\mu$ attains
its norm on $s_1$, so that $\Phi$ is an isometry.

For $\mu\in\mathbb M_m(X')$, by definition,
\begin{align*} \|\mu\| &=
\sup\{ |\langle\ip{\mu}{x}\rangle| : n\in\mathbb N,
x\in\mathbb M_n(X), \|x\|=1 \}
= \sup\{ |\langle\ip{\mu}{x}\rangle| : x\in s \}. \end{align*}
Following the notation in \cite{ER}, for $x=(x_{ij})\in\mathbb M_n(X)$
and $\mu=(\mu_{kl})\in\mathbb M_m(X')$, we let $\langle\ip{\mu}{x}\rangle
= ( \ip{\mu_{kl}}{x_{ij}} )_{(k,i),(l,j)} \in \mathbb M_m\otimes\mathbb M_n
= \mathbb M_{m\times n}$.
We then see that $(\Phi)_n(\mu) = ( \langle\ip{\mu}{x}\rangle )_{x\in s}$,
so that $(\Phi)_n$ is an isometry, and hence $\Phi$ is a $p$-complete
isometry as required.
\end{proof}

We now come to our first problem.  Let $X$ be a Banach space, and
recall the isometric map $\kappa = \kappa_X:X\rightarrow X''$ defined by
$\ip{\kappa_X(x)}{\mu} = \ip{\mu}{x}$ for $x\in X$ and $\mu\in X'$.

\begin{proposition}\label{kappa_iso}
Let $X$ be a $p$-operator space.  Then $\kappa_X$ is a $p$-complete
contraction.  Furthermore, $\kappa_X$ is a $p$-complete isometry
if and only if $X\subseteq\mc B(L_p(\phi))$ $p$-completely
isometrically for some measure $\phi$.
\end{proposition}
\begin{proof}
For $x=(x_{ij})\in\mathbb M_n(X)$, by definition,
\begin{align*}
\|(\kappa)_n(x)\|_n &= \sup\big\{ \|\langle\ip{(\kappa(x_{ij}))}{\mu}\rangle\| :
   m\in\mathbb N, \mu\in\mathbb M_m(X'), \|\mu\|_m=1 \big\} \\
&= \sup\big\{ \|\langle\ip{\mu}{x}\rangle\| :
   m\in\mathbb N, \mu\in\mathbb M_m(X'), \|\mu\|_m=1 \big\} \leq \|x\|_n,
\end{align*}
so that $\kappa$ is a $p$-complete contraction.

Suppose now that $\kappa$ is a $p$-complete isometry.  From the above
theorem, we know that $X'' \subseteq \mc B(\ell_p(I))$ for some index
set $I$.  Thus $X = \kappa(X) \subseteq \mc B(\ell_p(I))$, as required.

Conversely, suppose that $X\subseteq\mc B(E)$ for $E=L_p(\phi)$ for some
measure $\phi$.  To show that $\kappa$ is a $p$-complete isometry, we need to
show that for each $x\in\mathbb M_n(X)$ and $\epsilon>0$ there exists
$m\in\mathbb N$ and a $p$-complete contraction $u\in\mc{CB}_p(X,\mathbb M_m)
= \mathbb M_m(X')$ with $\| u(x) \|_m \geq \|x\|_n-\epsilon$.
When $p=2$, we may use \cite[Lemma~2.3.4]{ER} and take $n=m$ and
$\epsilon=0$.  However, for other values of $p$ we have to work harder.

Let $x=(x_{ij})\in\mathbb M_n(X)\subseteq\mc B(\ell_p^n\otimes_p E)$.  For
$\epsilon>0$, there exists $(a_i)_{i=1}^n\subseteq E$ with
$\sum_{i=1}^n \|a_i\|^p\leq 1$ and
\[ \Big(\sum_{i=1}^n \Big\| \sum_{j=1}^n x_{ij}(a_j) \Big\|^p \Big)^{1/p}
\geq \|x\|_n-\epsilon. \]
Let $b_i = \sum_{j=1}^n x_{ij}(a_j)$ for $1\leq i\leq n$.  Let
$\delta>0$ to be chosen later.
By standard properties of $E=L_p(\phi)$, there exists $m\in\mathbb N$
and an isometry $U:\ell_p^m \rightarrow E$ such that for each $j$,
there exists $f_j\in\ell_p^m$ with $\| U(f_j) - a_j \|<\delta$.
Similarly, there exists a contraction $V:E \rightarrow \ell_p^m$
such that $(1-\delta)\|b_i\| \leq \|V(b_i)\| \leq(1+\delta)\|b_i\|$
for each $i$.

Define $u:X\rightarrow\mc B(\ell_p^m)$ by $u(x) = VxU$ for $x\in X$.
A simple calculation shows that $u$ is a $p$-complete contraction, as
$\|U\|\|V\|\leq 1$.  Then
\[ \|(u)_n(x)\| \Big(\sum_{j=1}^n \|f_j\|^p \Big)^{1/p} \geq
\Big(\sum_{i=1}^n \Big\| \sum_{j=1}^n Vx_{ij}U(f_j) \Big\|^p \Big)^{1/p}
\geq \|x\|_n - 2\epsilon, \]
if $\delta>0$ is sufficiently small.  Similarly, if $\delta>0$ is
sufficiently small, then $\sum_{j=1}^n \|U(f_j)\|^p \leq \sum_{j=1}^n
\|a_j\|^p + \epsilon \leq 1+\epsilon$.  Hence $\|(u)_n(x)\|_m$ can be
chosen to be arbitrarily close to $\|x\|_n$, as required.
\end{proof}

The following was communicated to us by Christian Le Merdy.  Suppose
that $X\subseteq\mc B(L_p(\phi))$ for some measure $\phi$, and that $X$ is
finite dimensional with $\mathbb M_{n,1}(X) = \ell_p^n(X)$ for each $n$.
Pick $\epsilon>0$, and let $(x_1,\ldots,x_n)$ be an $\epsilon$-dense
subset of the unit sphere of $X$ (which exists as $X$ is finite
dimensional).  Then
\[ \Big( \sum_{k=1}^n \|x_k\|^p\Big)^{1/p}
= \| (x_k) \|_{\mathbb M_{n,1}(X)}
:= \sup\Big\{ \Big(\sum_{k=1}^n \|x_k(w)\|^p \Big)^{1/p} :
w\in L_p(\phi), \|w\|\leq1 \Big\}. \]
There hence exists $w_\epsilon\in L_p(\phi)$ with $\|w_\epsilon\|=1$ and
$\|x_k(w_\epsilon)\|\geq \|x_k\|-\epsilon=1-\epsilon$ for each $k$.  Define
$T_\epsilon:X \rightarrow L_p(\phi)$ by $T(x) = x(w_\epsilon)$ for $x\in X$.
For $x\in X$ with $\|x\|=1$, let $\|x-x_k\|<\epsilon$, so that
\[ 1=\|x\| \geq \|T_\epsilon(x)\| = \|x(w_\epsilon)\|
> \|x_k(w_\epsilon)\| - \epsilon \geq 1-2\epsilon. \]
By homogeneity, $(1-2\epsilon)\|x\| \leq \|T_\epsilon(x)\| \leq \|x\|$
for each $x\in X$.  A simple ultrapower argument then shows that
we may construct an isometry $X\rightarrow L_p(\psi)$ for some measure
$\psi$ (recall that an ultrapower of $L_p(\phi)$ is equal to
$L_p(\psi)$ for some $\psi$).

Now let $E \subseteq \ell_p^m$ be some subspace.  We give $\ell_p^m$
the $p$-operator space structure given by the identification
$\ell_p^m = \mc B(\mathbb C,\ell_p^m)$, and then make $E$ a subspace.
Then $\mathbb M_{n,1}(\ell_p^m) \subseteq \mc B(\mathbb C,\ell_p^m
\otimes_p \ell_p^n)$, so that $\mathbb M_{n,1}(\ell_p^m) =
\ell_p^n(\ell_p^m)$, and similarly for $E$.  In particular,
\[ \mathbb M_{n,1}(\ell_p^m / E) = \mathbb M_{n,1}(\ell_p^m)
/ \mathbb M_{n,1}(E) = \ell_p^n(\ell_p^m) / \ell_p^n(E)
= \ell_p^n(\ell_p^m/E). \]
So, if $\ell_p^m/E\subseteq\mc B(L_p(\phi))$ for some measure $\phi$,
then $\ell_p^m/E\subseteq L_p(\psi)$ for some measure $\psi$.  However,
for suitable chosen $E$, this is nonsense.  In particular, there
exist $p$-operator spaces $X$ (which may be finite-dimensional) such that
$\kappa_X$ is not a $p$-complete isometry.

\begin{lemma}
Let $X$ and $Y$ be $p$-operator spaces, and let $u\in\mc{CB}_p(X,Y)$.
Then $u'\in\mc{CB}_p(Y',X')$ and $\|u'\|_{pcb} \leq \|u\|_{pcb}$.
\end{lemma}
\begin{proof}
This follows as for operator spaces, see \cite[Proposition~3.2.2]{ER}.
We cannot conclude that $\|u'\|_{pcb} = \|u\|_{pcb}$ because of
the problems we encountered above.
\end{proof}

Combining Theorem~\ref{dual_embed} and Proposition~\ref{kappa_iso},
we see that for every $p$-operator space $X$, we have that
$\kappa_{X'}:X'\rightarrow X'''$ is a $p$-complete isometry.  Actually,
there is a much easier way to see this result.  A simple calculation
shows that $\kappa_X' \kappa_{X'} = I_{X'}$, and as the identity map
if a $p$-complete isometry, so also must $\kappa_{X'}$ be, as by the
lemma, $\kappa_X'$ is a $p$-complete contraction.

Let $X$ and $Y$ be $p$-operator spaces, and let $u\in\mc{CB}_p(X,Y)$.
The $u$ is a \emph{$p$-complete quotient map} if, for each $n$, $(u)_n$ takes
the open unit ball of $\mathbb M_n(X)$ onto the open unit ball of
$\mathbb M_n(Y)$.

\begin{lemma}\label{quot_to_iso}
Let $X$ and $Y$ be $p$-operator spaces, and let $u:X\rightarrow Y$ be a
$p$-complete quotient map.  Then $u':Y'\rightarrow X'$ is a $p$-complete
isometry.
\end{lemma}
\begin{proof}
Let $\mu\in\mathbb M_n(Y')$ and $\epsilon>0$, so that for some $m$,
there exists $y\in\mathbb M_m(Y)$ with $\|y\|_m<1$ and
$|\langle\ip{\mu}{y}\rangle|\geq \|\mu\|_n - \epsilon$.  By assumption,
we can find $x\in\mathbb M_m(X)$ with $\|x\|_m<1$ and $u(x)=y$, and so
\[ \|(u')_n(\mu)\|_n \geq \|\langle\ip{\mu}{u(x)}\rangle\|
\geq \|\mu\|_n-\epsilon, \]
which, as $\epsilon>0$ was arbitrary, shows that $\|(u')_n(\mu)\|_n = \|\mu\|_n$,
as required.
\end{proof}

The lack of a suitable Hahn-Banach theorem for $p$-operator spaces
(when $p=2$ we have the Arveson-Wittstock theorem \cite[Theorem~4.1.5]{ER})
means that we cannot show the converse to the above.

We define subspaces of $p$-operator spaces in the obvious way.  Given
a $p$-operator space $X$ and a closed subspace $Y\subseteq X$, we define
a norm on $\mathbb M_n(X/Y)$ by identifying this space with
$\mathbb M_n(X) / \mathbb M_n(Y)$.  Then, as for operator spaces
(see \cite[Proposition~3.11]{ER}) it is easy to check that $X/Y$ becomes
a $p$-operator space, and that the quotient map $\pi:X\rightarrow X/Y$ is
a $p$-complete quotient map.  The above lemma then tells us that
$\pi':(X/Y)'\rightarrow X'$ is a $p$-complete isometry.  A simple
calculation shows that the image of $\pi'$ is
\[ Y^\perp := \{ \mu\in X' : \ip{\mu}{y}=0 \ (y\in Y) \}, \]
so that we may identify $(X/Y)'$ with $Y^\perp$ $p$-completely isometrically.
Again, we have no such identification of $Y'$ with a suitable
quotient of $X'$.

\subsection{Tensor products}

We define the \emph{$p$-operator space projective tensor norm}
on the tensor product of two $p$-operator space $X$ and $Y$ to be
\[ \|\tau\|_{\wedge} = \inf\Big\{ \|\alpha\| \|u\| \|v\| \|\beta\|
: \tau = \alpha(u\otimes v)\beta \Big\} \qquad ( \tau\in\mathbb M_n(X\otimes Y) ). \]
Here we let $u\in\mathbb M_r(X)$ and $v\in\mathbb M_s(Y)$, so that
$u\otimes v \in\mathbb M_{r\times s}(X\otimes Y)$ in a natural way, and
we take $\alpha\in\mathbb M_{n,r\times s}$ and $\beta\in
\mathbb M_{r\times s,n}$, so that $\alpha(u\otimes v)\beta\in
\mathbb M_n(X\otimes Y)$ as required.  This is exactly the definition
for operator spaces, except that as above, we evaluate $\|\alpha\|$ as
a member of $\mc B(\ell_p^n,\ell_p^{s\times r})$, and similarly
$\|\beta\|$.  We shall prove below
that $\|\cdot\|_{\wedge}$ gives $X\otimes Y$ an abstract $p$-operator space
structure.  Denote by $X \proten^p Y$ the completion.

\begin{proposition}
Let $X$ be a vector space, and for each $n$, let $\|\cdot\|_n:
\mathbb M_n(X) \rightarrow [0,\infty)$ be a map such that:
\begin{itemize}
\item[$\mc D_\infty'$] for $u\in\mathbb M_n(X)$ and $v\in\mathbb M_m(X)$,
   we have that $\|u\oplus v\|_{n+m} \leq \max(\|u\|_n, \|v\|_m)$;
\item[\Axtwo] for $u\in\mathbb M_m(X)$, $\alpha\in\mathbb M_{n,m}$
   and $\beta\in\mathbb M_{m,n}$, we have that $\|\alpha u \beta\|_n
   \leq \|\alpha\| \|u\|_m \|\beta\|$.
\end{itemize}
Then each $\|\cdot\|_n$ is a norm, and the completion of $X$ becomes
an abstract $p$-operator space.
\end{proposition}
\begin{proof}
This follows exactly as for operator spaces, \cite[Proposition~2.3.6]{ER}.
\end{proof}

\begin{proposition}
Let $X$ and $Y$ be $p$-operator spaces.  Then $\|\cdot\|_{\wedge}$
induces a $p$-operator space structure on $X\otimes Y$.  Furthermore,
$\|\cdot\|_{\wedge}$ is the largest such $p$-operator space norm with
the additional property that $\|u\otimes v\| \leq \|u\|_r \|v\|_s$ for
$u\in\mathbb M_r(X)$ and $v\in\mathbb M_s(Y)$.
\end{proposition}
\begin{proof}
This follows as for operator space (see \cite[Theorem~7.1.1]{ER}) with minor
alterations.  In \cite{ER}, the authors use the C$^*$-identity, in the $p=2$
case, to estimate the norm of a matrix $\alpha\in\mathbb M_{r,s} =
\mc B(\ell_p^s,\ell_p^r)$ of the block form
\[ \alpha = \begin{pmatrix} \alpha_1 & 0 & 0 & 0 \\ 0 & 0 & 0 & \alpha_2
\end{pmatrix}. \]
However, we get the entirely elementary estimate that $\|\alpha\|
\leq \max ( \|\alpha_1\|, \|\alpha_2\|)$, which is all that is
required.
\end{proof}

Let $X, Y$ and $Z$ be $p$-operator spaces, and let $\psi:X\times Y
\rightarrow Z$ be a bilinear map.  We define bilinear maps
\[ (\psi)_{r,s;t,u} : \mathbb M_{r,s}(X)\times\mathbb M_{t,u}(Y)
\rightarrow \mathbb M_{r\times t,s\times u}(Z); \ \
(x,y) \mapsto \big( \psi(x_{i,j},y_{k,l}) \big). \]
Then we let $(\psi)_{r;s} = (\psi)_{r,r;s,s}$, and define
\[ \|\psi\|_{pcb} = \sup\{ \|(\psi)_{r;s}\| : r,s\in\mathbb N \}. \]
This leads to the definition of $\mc{CB}_p(X\times Y,Z)$, which
can be turned into a $p$-operator space in the same way as for $\mc{CB}_p$.

\begin{proposition}
Let $X,Y$ and $Z$ be operator spaces.  Then we have natural completely
isometric identifications
\[ \mc{CB}_p(X\proten^p Y,Z) = \mc{CB}_p(X\times Y,Z)
= \mc{CB}_p(X,\mc{CB}_p(Y,Z)). \]
\end{proposition}
\begin{proof}
This follows as for operator spaces, see \cite[Proposition~7.1.2]{ER}.
\end{proof}

We hence see that, for example, $(X\proten^p Y)' = \mc{CB}_p(X,Y')$.
As for operator spaces (see \cite[Chapter~7]{ER}), we can now easily
show that $X\proten^p Y = Y\proten^p X$ naturally, and that the
operator $\proten^p$ is associative.  Furthermore, if $u_i:X_i\rightarrow
Y_i$ are complete contractions for $i=1,2$, then $u_1\otimes u_2$
extends to a complete contraction $X_1\proten^p X_2 \rightarrow
Y_1\proten^p Y_2$.

\begin{proposition}\label{proj_is_proj}
Let $X,Y,X_1$ and $Y_1$ be $p$-operator spaces, and let
$u:X\rightarrow X_1$ and $v:Y\rightarrow Y_1$ be $p$-complete quotient
maps.  Then
$u\otimes v:X\proten^p Y \rightarrow X_1\proten^p Y_1$ is also a
$p$-complete quotient map.  Furthermore, $\ker(u\otimes v)$ is the
closure of the space
\[ (\ker u) \otimes Y + X \otimes (\ker v) \subseteq X\proten^p Y. \]
\end{proposition}
\begin{proof}
A careful examination of the proof for operator spaces,
\cite[Proposition~7.1.7]{ER}, shows that the proof is equally valid
for $p$-operator spaces.
\end{proof}

\section{Algebras}

In this section, we shall study weak$^*$-closed subalgebras of
$\mc B(E)$ for an $SQ_p$ space $E$.  The starting point is to look
at $\mc B(E)$ itself, and in particular, its predual $E'\proten E$.

Let $\phi$ be a measure, and consider the space $\mc N(L_p(\phi))$
of nuclear operators on $L_p(\phi)$, so that $\mc N(L_p(\phi))' =
\mc B(L_p(\phi))$ as explained above.  Thus $\mc N(L_p(\phi))$ carries
a natural $p$-operator space structure by duality.

\begin{lemma}\label{nuc_dual_op}
With notation as above, $\mc B(L_p(\phi)) = \mc N(L_p(\phi))'$
$p$-completely isometrically.
\end{lemma}
\begin{proof}
To ease notation, write $\mc N = \mc N(L_p(\phi))$ and
$\mc B = \mc B(L_p(\phi))$.  By definition, for $\tau\in\mathbb M_n(\mc N)$,
we have that
\[ \|\tau\|_n = \sup\{\|\langle\ip{\tau}{T}\rangle\| : m\in\mathbb N,
T\in\mathbb M_m(\mc B), \|T\|_m\leq1\}. \]
Here we have identified $\mathbb M_n(\mc N)$ with a subspace of
$\mathbb M_n(\mc B') = \mc{CB}_p(\mc B,\mathbb M_n)$, and it is easy to
see that this subspace coincides with the space $\mc{CB}_p^\sigma(\mc B,
\mathbb M_n)$ of weak$^*$-continuous $p$-completely bounded maps from
$\mc B$ to $\mathbb M_n$.

For $T\in\mathbb M_n(\mc B)$, let $\|T\|_{\mc N'}$ be the norm
of $T$ considered as a member of $\mathbb M_n(\mc N') =
\mc{CB}_p(\mc N,\mathbb M_n)$, so that
\[ \|T\|_{\mc N'} = \sup\{\|\langle\ip{\tau}{T}\rangle\| : m\in\mathbb N,
\tau\in\mathbb M_m(\mc N), \|\tau\|\leq1\} \leq \|T\|. \]
To show the converse, for $\epsilon>0$, we wish to find $\tau\in
\mathbb M_m(\mc N) = \mc{CB}_p^\sigma(\mc B,\mathbb M_n)$ with
$|\langle\ip{\tau}{T}\rangle| \geq \|T\|-\epsilon$.

By Proposition~\ref{kappa_iso}, we know that there exists
$\tau\in\mc{CB}_p(\mc B,\mathbb M_m)$ with this property.  Following that
proof, we see that $\tau$ is defined to be $\tau(T) = VTU$ for $T\in\mc B$,
for suitable $U:\ell_p^m\rightarrow L_p(\phi)$ and $V:L_p(\phi)\rightarrow\ell_p^m$.
A simple calculation shows that such a map is actually in
$\mc{CB}_p^\sigma(\mc B,\mathbb M_n)$, which completes the proof.
\end{proof}

It will be useful to have a more concrete description of the norm
on $\mc N(L_p(\phi))$.  For ease of notation, let $\mc N=\mc N(L_p(\phi))$
and $\mc B=\mc B(L_p(\phi))$.  Let $n\in\mathbb N$ and
$\tau\in\mathbb M_n(\mc N)$.  Then, as above, $\|\tau\|_n = 
\sup\{ \|\langle\ip{T}{\tau}\rangle\| : T\in\mathbb M_m(\mc B),
\|T\|_m\leq 1\}$.  For $T\in\mathbb M_m(\mc B)$, we have that
\[ \|\langle\ip{T}{\tau}\rangle\| = \sup\Big\{ \Big|
\sum_{i=1}^n \sum_{k=1}^m \sum_{j=1}^n \sum_{l=1}^m \beta_{ki}
\ip{T_{kl}}{\tau_{ij}} \alpha_{lj} \Big|
: \sum_{l,j} |\alpha_{lj}|^p\leq 1, \sum_{k,i} |\beta_{ki}|^{p'}\leq1 \Big\}. \]
Suppose that $\tau_{ij} = \sum_{r=1}^\infty \mu_r^{(ij)}
\otimes x_r^{(ij)} \in L_{p'}(\phi) \proten L_p(\phi)$ for each $i,j$.
Treat $T = (T_{kl})\in\mathbb M_m(\mc B)$ as an operator on $\ell_p^m\otimes_p
L_p(\phi)$, given by $T(\delta_l\otimes x) = \sum_{k=1}^m
\delta_k\otimes T_{kl}(x)$.  Then
\begin{align*}
\|\langle\ip{T}{\tau}\rangle\| &= \sup\Big\{ \Big| \sum_{i,j,k,l}
   \sum_{r=1}^\infty \beta_{ki} \ip{\mu_r^{(ij)}}{T_{kl}(x_r^{(ij)})} \alpha_{lj} \Big|
   : \sum_{l,j} |\alpha_{lj}|^p\leq 1, \sum_{k,i} |\beta_{ki}|^{p'}\leq1 \Big\} \\
&\hspace{-5ex}= \sup\Big\{ \Big| \sum_{i,j,k,l} \sum_{r=1}^\infty
   \ip{\beta_{ki}\delta_k^*\otimes\mu_r^{(ij)}}{T(\alpha_{lj}\delta_l\otimes x_r^{(ij)})}
   \Big| : \sum_{l,j} |\alpha_{lj}|^p\leq 1, \sum_{k,i} |\beta_{ki}|^{p'}\leq1 \Big\} \\
&\hspace{-5ex}= \sup\Big\{ \Big| \sum_{i,j} \sum_{r=1}^\infty
   \ip{\eta_i\otimes\mu_r^{(ij)}}{T(\gamma_j\otimes x_r^{(ij)})}
   \Big| : \sum_j \|\gamma_j\|^p\leq1, \sum_i \|\eta_i\|^{p'}\leq1 \Big\},
\end{align*}
where we have $(\eta_i)\subseteq\ell_{p'}^m$ and $(\gamma_j)\subseteq\ell_p^m$.
Thus, by the usual duality between $\mc N(\ell_p^m\otimes_p L_p(\phi))$ and
$\mc B(\ell_p^m\otimes_p L_p(\phi))$, we see that
\begin{equation}
\|\tau\|_n = \sup\Big\{ \Big\| \sum_{r=1}^\infty \sum_{i,j=1}^n
   \big( \eta_i\otimes\mu_r^{(ij)} \big) \otimes
   \big( \gamma_j\otimes x_r^{(ij)} \big) \Big\|_\pi
: \sum_i \|\eta_i\|^{p'}\leq1, \sum_j \|\gamma_j\|^p\leq1 \Big\},
\label{eq:one} \end{equation}
where now $m$ is also free to vary.

Let $\mc N_n^p = \mc N(\ell_p^n)$, so by the lemma, $(\mc N_n^p)' =
\mc B(\ell_p^n) = \mathbb M_n$.  For
a $p$-operator space $X$, we hence have that
\[ (\mc N_n^p \proten^p X)' = \mc{CB}_p(X,(\mc N_n^p)')
= \mc{CB}_p(X,\mathbb M_n) = \mathbb M_n(X'). \]
In particular,
\[ (\mc N_n^p \proten^p \mc N_m^p)' = \mathbb M_n((\mc N_m^p)')
= \mathbb M_n(\mathbb M_m) = \mathbb M_{n\times m}, \]
and so, as everything is finite-dimensional,
\[ \mc N_n^p \proten^p \mc N_m^p = \mc N_{n\times m}^p, \]
completely isometrically.

\begin{proposition}\label{nuc_p_proj_prod}
We have a natural completely isometric identification
\[ \mc N(\ell_p) \proten^p \mc N(\ell_p) = \mc N(\ell_p \otimes_p \ell_p). \]
\end{proposition}
\begin{proof}
We follow the proof of \cite[Proposition~7.2.1]{ER}.
For $n\in\mathbb N$, let $\iota_n:\ell_p^n\rightarrow\ell_p$ be
the inclusion onto the first $n$ co-ordinates, and let
$p_n:\ell_p\rightarrow\ell_p^n$ be the natural projection.  Thus
the maps
\begin{gather*} j_n:\mc N(\ell_p^n) \rightarrow \mc N(\ell_p);\
\tau\mapsto \iota_n \tau p_n, \\
P_n:\mc N(\ell_p)\rightarrow\mc N(\ell_p^n);\
\sigma\mapsto p_n \sigma \iota_n,
\end{gather*}
are, respectively, a complete isometry and a complete quotient map
such that $P_n j_n$ is the identity.  Thus $j_n P_n$ is a completely
contractive projection of $\mc N(\ell_p)$ onto $\mc N_n^p$.

For $n,m$, we have the commutative diagram
\[ \xymatrix{ \mc N_p^n \proten^p \mc N_p^m \ar[r]
\ar[d]^{j_n\otimes j_m} &
\mc N(\ell_p^n \otimes_p \ell_p^m) \ar[d] \\
\mc N(\ell_p) \proten^p \mc N(\ell_p) \ar[r] &
\mc N(\ell_p\otimes_p\ell_p) } \]
As above, we know that the top row is a complete isometry.  From
the previous paragraph, we know that $j_n\otimes j_m$ is a
complete isometry, and similarly, the right column is a complete
isometry.  The union of the spaces $\mc N^n_p\otimes\mc N^m_p$
is norm dense in $\mc N(\ell_p)\proten^p\mc N(\ell_p)$, and
the union of the spaces $\mc N(\ell^n_p\otimes_p\ell^m_p)$ is
norm dense in $\mc N(\ell_p\otimes_p\ell_p)$.  Hence, as all the
maps are coherent, we conclude that the bottom row must also
be a complete isometry, as required.
\end{proof}

\begin{proposition}\label{proten_nuc}
Let $\phi$ and $\lambda$ be measures.  We have a natural
completely isometric identification
\[ \mc N(L_p(\phi)) \proten^p \mc N(L_p(\lambda)) =
\mc N(L_p(\phi\times\lambda)). \]
\end{proposition}
\begin{proof}
Spaces of the form $L_p(\mu)$ admit a  net of subspaces $(E_i)$
whose union is dense, and such that each $E_i$ is $1$-complemented,
and isometric to $\ell_p^n$ for some $n$.  Hence we may directly
adapt the above proof.
\end{proof}

Suppose that such a net of subspaces $(E_i)$ exists for some $E\in SQ_p$.
Then it is easily seen that $E$ is a $\mc{L}_{p,1}^g$ space,
as defined in \cite[Section~23]{DF}.  By \cite[Theorem~23.2]{DF}, $E$ is
thus isometric to a $1$-complemented subspace of some $L_p$ space, and is
thus isometric to an $L_p$ space (see \cite{T}).  Hence the above proposition
is the best we can do, at least using this method of proof.

We wish to further study the norm on $\mathbb M_n(\mc N(E))$,
for $E\in SQ_p$.  Suppose that $E$ has the approximation property
(eventually, we shall have to assume that $E=L_p(\phi)$ anyway) so that
$\mc K(E)' = \mc N(E)$.  Define $T_n(\mc K(E))$ to be the vector space
$\mathbb M_n(\mc K(E))$ together with the norm defined by,
for $K=(k_{ij})_{i,j=1}^n$,
\[ \|K\|_{T_n(\mc K(E))} =
\inf\Big\{ \|T\|_m \Big(\sum_{i,k} |\alpha_{ik}|^p \Big)^{1/p}
\Big(\sum_{i,k} |\beta_{ik}|^{p'} \Big)^{1/p'} \Big\}, \]
where we take the infimum over $m\in\mathbb N$ and
$T\in\mathbb M_m(\mc K(E))$ such that for each $i,j$,
$k_{ij} = \sum_{k,l=1}^m \beta_{ki} T_{kl} \alpha_{lj}$.
We define a bilinear mapping $\mathbb M_n(\mc N(E)) \times
T_n(\mc K(E)) \rightarrow\mathbb C$ by
\[ \ip{\tau}{K} = \sum_{i,j=1}^n \ip{\tau_{ij}}{K_{ij}}
\qquad \big( \tau=(\tau_{i,j}) \in\mathbb M_n(\mc N(E)),
K=(k_{ij})\in T_n(\mc K(E)) \big). \]
By formula (\ref{eq:one}) it is immediate that $|\ip{\tau}{K}|
\leq \|\tau\|_n \|K\|_{T_n(\mc K(E))}$.

Let $\Gamma\in T_n(\mc K(E))'$, and for each $i,j$, define
$\tau_{ij}\in\mc N(E)$ by $\ip{\tau_{ij}}{k} = \ip{\Gamma}{
\delta_{ij}\otimes k}$ for $k\in\mc K(E)$.  Here $\delta_{ij}
\otimes k\in T_n(\mc K(E))$ is the matrix with $k$ in the $(i,j)$
entry, and $0$ elsewhere.  Then $\|\delta_{ij}\otimes k\|_{T_n(\mc K(E))}
\leq \|k\|$, so that $\tau_{ij}$ is well-defined, and
$\| \tau_{ij} \| \leq \|\Gamma\|$.  Let $\tau = (\tau_{ij})
\in \mathbb M_n(\mc N(E))$.  Let $\tau_{ij} = \sum_r
\mu_r^{(ij)} \otimes x_r^{(ij)}$ for each $i,j$.
Then let $T\in\mathbb M_m(\mc K(E))$, so that
\begin{align*}
\|\langle\ip{T}{\tau}\rangle\| &= \sup\Big\{ \Big| \sum_{i,j,k,l}
   \sum_{r=1}^\infty \beta_{ki} \ip{\mu_r^{(ij)}}{T_{kl}(x_r^{(ij)})} \alpha_{lj} \Big|
   : \sum_{l,j} |\alpha_{lj}|^p\leq 1, \sum_{k,i} |\beta_{ki}|^{p'}\leq1 \Big\} \\
&= |\ip{\tau}{K}|,
\end{align*}
where $K=(k_{ij})\in T_n(\mc K(E))$ is defined by
$k_{ij} = \sum_{k,l=1}^m \beta_{ki} T_{kl} \alpha_{lj}$.
By definition, $\|K\|_{T_n(\mc K(E))} \leq 1$, so by the definition
of $\mathbb M_n(\mc N(E))$, we conclude that $T_n(\mc K(E))' =
\mathbb M_n(\mc N(E))$ isometrically.  Here we move from taking
a supremum over $\mathbb M_m(\mc B(E))$ to $\mathbb M_m(\mc K(E))$,
which we may do by approximation, as $E$ has the (metric) approximation property.

Define $T_n(\mc B(E))$ in a similar way to the definition of
$T_n(\mc K(E))$.  Given $T=(T_{ij}) \in \mathbb M_n(\mc B(E))$
and $\tau=(\tau_{ij})\in\mathbb M_n(\mc N(E))$, so that we see that
$|\ip{T}{\tau}| \leq \|T\|_n \|\tau\|_n$ immediately.
Proceeding as above, we may at least identify $\mathbb M_n(\mc N(E))'$
with $T_n(\mc B(E))$ as vector spaces.

\begin{proposition}\label{dual_of_nuc}
Let $\phi$ be a measure, and let $E=L_p(\phi)$.
Then $\mathbb M_n(\mc N(E))' = T_n(\mc B(E))$ isometrically.
\end{proposition}
\begin{proof}
Suppose firstly that $E$ is finite-dimensional (that is,
$E=\ell_p^N$ for some $N$).  Then $\mc B(E) = \mc K(E)$,
and as the space $\mathbb M_n(\mc N(E))$ is finite-dimensional,
we see that $\mathbb M_n(\mc N(E))' = T_n(\mc B(E))$.
The general case then follows by a finite-dimensional
decomposition argument, as used in Proposition~\ref{nuc_p_proj_prod}.

Indeed, let $F\subseteq E$ be a $1$-complemented finite-dimensional
subspace.  Thus $\mc N(\ell_p^m \otimes_p F) \subseteq
\mc N(\ell_p^m \otimes_p E)$ isometrically, for each $m$.
It hence follows that $\mathbb M_n(\mc N(F)) \subseteq
\mathbb M_n(\mc N(E))$ isometrically, and so the natural map
$\mathbb M_n(\mc N(E))' \rightarrow \mathbb M_n(\mc N(F))'$ is a
quotient map.  Similarly, we may check that the natural map
$T_n(\mc B(E)) \rightarrow T_n(\mc B(F))$ (induced by the projection of
$E$ onto $F$) is a quotient map.  Thus we have the following diagram
\[ \xymatrix{ \mathbb M_n(\mc N(F))' &
\mathbb M_n(\mc N(E))' \ar[l] \\
T_n(\mc B(F)) \ar[u]^{\cong} & T_n(\mc B(E)). \ar[l]^{\phi_F} \ar[u]^\psi } \]
The map on the left is norm-decreasing, while the map on the right
is an isometric isomorphism.  Let $T\in T_n(\mc B(E))$, and we may
easily check that
\[ \|T\|_{T_n(\mc B(E))} = \sup \{ \| \phi_F(T) \|_{T_n(\mc B(F))}
: F\subseteq E \}. \]
The supremum is taken over $1$-complemented subspaces of $E$, of course.
A similar equality holds for $\psi(T)$, and hence it follows that
$\|\psi(T)\|_{\mathbb M_n(\mc N(E))'} = \|T\|_{T_n(\mc B(E))}$, as required.
\end{proof}

As before, this method of proof does not readily generalise to spaces
other than $L_p(\phi)$.

\subsection{General weak$^*$-closed algebras}\label{weakstaralgs}

Let $E=L_p(\phi)$ for some measure $\phi$, and let $\mc A\subseteq
\mc B(E)$ be a weak$^*$-closed algebra.  The predual of $\mc A$,
denoted $\mc A_*$, may be identified with the quotient
$\mc A_* = \mc N(E) / {^\perp}\mc A$, where
\[ {^\perp}\mc A = \{ \tau\in\mc N(E) : 
\ip{a}{\tau}=0 \ (a\in\mc A) \}. \]
Clearly $\mc A$ carries a canonical $p$-operator space structure,
and we can use this to induce a $p$-operator space structure on $\mc A_*$.
We shall call this the \emph{dual structure} on $\mc A_*$.

\begin{proposition}\label{com_duality_alg}
Let $\mc A\subseteq\mc B(L_p(\phi))$ be a weak$^*$-closed subalgebra,
for some measure $\phi$.  Give $\mc A_*$ the dual structure.
Then $\mc A_*' = \mc A$ $p$-completely isometrically.
\end{proposition}
\begin{proof}
This follows in an analogous way to the proof of Lemma~\ref{nuc_dual_op}.
To be precise, let $T\in\mathbb M_n(\mc A)$ and $\epsilon>0$.
Then there exists $m\in\mathbb N$ and maps $U:\ell_p^m\rightarrow L_p(\phi)$
and $V:L_p(\phi)\rightarrow\ell_p^m$ such that $\|U\|=\|V\|=1$ and,
if $\tau\in\mc{CB}_p(\mc A,\mathbb M_m)$ is defined by $\tau(a) = VaU$,
then $\|\langle\ip{a}{\tau}\rangle\| \geq (\|a\|-\epsilon)\|\tau\|$.

Define $\sigma\in\mathbb M_m(\mc A_*) = \mc{CB}_p^\sigma(\mc A,
\mathbb M_m)$ by setting
\[ \sigma_{ij} = \tau_{ij} + {^\perp}\mc A \in
\mc N(L_p(\phi)) / {^\perp}\mc A = \mc A_*
\qquad (1\leq i,j\leq m). \]
Then $\langle\ip{a}{\tau}\rangle = \langle\ip{a}{\sigma}\rangle$,
and we claim that $\|\sigma\| \leq \|\tau\|$, which will complete the
proof.  To show this claim, it suffices to show that as an operator
in $\mc{CB}_p(\mc A,\mathbb M_m)$, $\sigma$ is a contraction.  This
is immediate however, as $\sigma$ agrees with $\tau$ on $\mc A$.
\end{proof}

Notice $\mc A_*$ is also a quotient of $\mc N(E)$, and so we could
define a $p$-operator space structure on $\mc A_*$ by insisting that
the quotient map $\pi:\mc N(E)\rightarrow\mc A_*$ is a $p$-complete
quotient mapping.  We shall call this the \emph{quotient structure}.
By Lemma~\ref{quot_to_iso}, when $\mc A_*$ has the quotient structure,
the inclusion $\pi':\mc A = \mc A_*' \rightarrow \mc N(E)'=\mc B(E)$ is
a $p$-complete isometry.  Thus $\mc A$ carries the same $p$-operator
space structure, irrespective of the $p$-operator space structure
put on $\mc A_*$.  We also see that, in general, the quotient norm
dominates the dual norm on $\mathbb M_n(\mc A_*)$ for each $n$.
When $p=2$, we may immediate conclude that
the two structures on $\mc A_*$ coincide, but for other values of $p$,
the lack of a suitable Hahn-Banach result means that we cannot conclude
this.  We shall later show that this problem seems to have some
link with amenability (see Theorem~\ref{when_amen_thm}), a result we
prepare for now.

Let $E=L_p(\phi)$ for some measure $\phi$.
From Proposition~\ref{dual_of_nuc}, we know that
$\mathbb M_n(\mc N(E))' = T_n(\mc B(E))$ isometrically.  We may regard
$(\pi)_n$ as a map from $\mathbb M_n(\mc N(E))$ to $\mathbb M_n(\mc A_*)$,
which is defined to be a quotient map when $\mc A_*$ carries the quotient
structure.  Thus $(\pi)_n' : \mathbb M_n(\mc A_*)'\rightarrow
\mathbb M_n(\mc N(E))' = T_n(\mc B(E))$ is an isometry which maps onto
$( \ker(\pi)_n)^\perp$.  It is easy to see that $\tau\in\mathbb M_n(\mc N(E))$
lies in $\ker(\pi)_n$ if and only if $\tau_{ij}\in\ker\pi$ for
each $i,j$.  Hence it follows that $T\in T_n(\mc B(E))$ lies in the image
of $(\pi)_n'$ if and only if $T_{ij}\in \mc A$ for each $i,j$.

From the definition of $T_n(\mc B(E))$, we see that the quotient structure
norm on $\mathbb M_n(\mc A_*)$ may be computed by considering matricies
$T=(T_{ij})$ such that $T_{ij} = \sum_{k,l} \beta_{ik} S_{kl} \alpha_{lj}
\in \mc A$ for some $S\in\mathbb M_m(\mc B(E))$ of norm one, and suitable
$\alpha$ and $\beta$.  By definition, the dual structure norm may be computed
by exactly the same method, only now we must ensure that $S_{kl}\in\mc A$
for each $k,l$, and not only that $T_{ij}\in\mc A$ for each $i,j$.

\begin{proposition}\label{proj_implies_struc_equal}
Let $\mc A$ and $\mc A_*$ be as above, and suppose that there is a
$p$-completely contractive projection from $\mc B(E)$ onto $\mc A$.
Then the two $p$-operator space structures on $\mc A_*$ coincide.
\end{proposition}
\begin{proof}
This is immediate, as given $T=(T_{ij})$ with $T_{ij} = \sum_{k,l}
\beta_{ik} S_{kl} \alpha_{lj} \in \mc A$ for each $i,j$, then we have
that $P(T_{ij}) = \sum_{k,l} \beta_{ik} P(S_{kl}) \alpha_{lj} \in \mc A$,
where $P(S_{kl})\in\mc A$ for each $k,l$.  As $\| (P(S_{kl})) \|_n
= \| (P)_n (S) \|_n \leq \|P\|_{pcb} \|S\|_n = \|S\|_n$, the claim follows.
\end{proof}

\section{Tensor products of algebras}\label{algebras}

For two von Neumann algebras $R$ and $S$, there is a natural
tensor product of their preduals $R_*$ and $S_*$ such that
$R_* \otimes S_*$ is the predual of the von Neumann algebra tensor
product $R\vnten S$.  A key fact about operator spaces
(\cite[Theorem~7.2.4]{ER}) is that $R_* \proten^2 S_*$ agrees with
the predual of $R\vnten S$.  In this section, we shall explore how
this result is proved, and shall lay the foundations for analogous
proofs, in the $p\not=2$ case, in some rather special cases.

We shall now study Slice Maps, following the presentation in
\cite[Section~7.2]{ER}.  Let $\phi_1,\phi_2$ be measures, and set
$E=L_p(\phi_1)$ and $F=L_p(\phi_2)$.
Let $w_1\in\mc N(E)$, so that we have a map $w_1\otimes I:
\mc B(E)\otimes\mc B(F)\rightarrow\mc B(F)$ given by
$(w_1\otimes I)(T\otimes S) = \ip{T}{w_1}S$.

\begin{lemma}\label{slice_map_lemma}
There exists a weak$^*$-continuous map $R(w_1):\mc B(E\otimes_p F)
\rightarrow\mc B(F)$ such that $R(w_1)$, when restricted to
$\mc B(E)\otimes\mc B(F)$, agrees with $w_1\otimes I$.  Furthermore,
$R(w_1)$ is $p$-completely bounded with $\|R(w_1)\|_{pcb} = \|w_1\|$.
\end{lemma}
\begin{proof}
For $u\in\mc B(E\otimes_p F)$, define 
$R(w_1)(u)\in\mc B(F) = \mc N(F)'$ by
\[ \ip{R(w_1)(u)}{\tau} = \ip{u}{w_1\otimes\tau}
\qquad (\tau\in\mc N(F)). \]
Then clearly $R(w_1)(u)\in\mc B(F)$ and $\|R(w_1)(u)\| \leq
\|u\| \|w_1\|$.  Obviously $R(w_1):\mc B(E\otimes_p F)\rightarrow
\mc B(F)$ is linear, and is thus a
bounded operator which clearly extends $w_1\otimes I$.
Furthermore, we may define $r(w_1):\mc N(F)\rightarrow\mc N(E\otimes_p F)$
by
\[ r(w_1)(\tau) = w_1\otimes\tau \in \mc N(E)\otimes\mc N(F)
\subseteq\mc N(E\otimes_p F) \qquad (\tau\in\mc N(F)), \]
and then we clearly see that $r(w_1)' = R(w_1)$, so that $R(w_1)$
is weak$^*$-continuous.

By Proposition~\ref{proten_nuc}, $\mc N(E)\proten^p\mc N(F) = 
\mc N(E\otimes_p F)$, and so $\mc B(E\otimes_p F) = 
\mc{CB}_p(\mc N(E),\mc B(F))$ $p$-completely isometrically.
Concretely, this second identification is given as follows.
For $u\in\mc B(E\otimes_p F)$, we define $\Lambda(u)\in
\mc{CB}_p(\mc N(E),\mc B(F))$ by
\[ \Lambda(u)(w_1) = R(w_1)(u) \qquad (w_1\in\mc N(E)). \]
Let $U \in \mathbb M_n(\mc B(E\otimes_p F))$ so that
$(R(w_1))_n(U) \in \mathbb M_n(\mc B(F))$.  Then
\[ (R(w_1))_n(U) = ( R(w_1)(U_{ij}) ) = ( \Lambda(U_{ij})(w_1) )
= (\Lambda)_n(U)(w_1), \]
so that $\|(R(w_1))_n(U)\| = \|(\Lambda)_n(U)(w_1)\|
\leq \|\Lambda\|_{pcb} \|U\| \|w_1\| = \|U\| \|w_1\|$, and
so $\|(R(w_1))_n\| \leq \|w_1\|$, implying that $\|R(w_1)\|_{pcb}
\leq \|w_1\|$.  Clearly then $\|R(w_1)\|_{pcb} = \|w_1\|$, as required.
\end{proof}

Similarly, we may work ``on the left'', leading to the definition
of $L(w_2):\mc B(E\otimes_p F) \rightarrow \mc B(E)$ for $w_2\in
\mc N(F)$.

Given weak$^*$-closed subalgebras
$\mc A\subseteq\mc B(E)$ and $\mc B\subseteq\mc B(F)$, we define
$\mc A\vnten\mc B$ to be the weak$^*$-closure of $\mc A\otimes\mc B$
in $\mc B(E\otimes_p F) = \mc B(L_p(\phi_1\times\phi_2))$.  We define
the \emph{Fubini product} $\mc A \otimes_{\mc F} \mc B$ to be the
subspace
\[ \{ u\in\mc B(E\otimes_p F) : R(w_1)(u)\in\mc B, L(w_2)(u)\in\mc A
\ (w_1\in\mc N(E), w_2\in\mc N(F)) \}. \]
As $R(w_1)$ and $L(w_2)$ are weak$^*$-continuous, we immediately see
that $\mc A\vnten\mc B \subseteq\mc A \otimes_{\mc F}\mc B$.

In general, we can only say a little about $\mc A\vnten\mc B$.
Let $w_1\in\mc N(E)$, and consider the map $R(w_1)$ restricted to
$\mc A\vnten\mc B$, which by weak$^*$-continuity maps into $\mc B$.
Suppose that $w_2\in\mc N(E)$ is such that $w_1-w_2\in{^\perp}\mc A$.
Then, for any $\tau\in\mc N(F)$, clearly $(w_1-w_2)\otimes\tau$
annihilates $\mc A\otimes\mc B$, and so
\[ \ip{R(w_1-w_2)(T)}{\tau} = \ip{T}{(w_1-w_2)\otimes\tau} = 0
\qquad (T\in\mc A\vnten \mc B). \]
Hence $R$ becomes a well-defined map $\mc N(E)/{^\perp}\mc A = \mc A_*
\rightarrow\mc{CB}_p(\mc A\vnten\mc B,\mc B)$, and similarly for $L$.

Now define a map $\delta:\mc A\vnten\mc B\rightarrow (\mc A_* \proten^p \mc B_*)'
= \mc{CB}_p(\mc B_*,\mc A)$ by
\[ \ip{\delta(T)}{\tau\otimes\sigma} =
\ip{R(\tau)(T)}{\sigma} = \ip{L(\sigma)(T)}{\tau}
\qquad (T\in\mc A\vnten\mc B, \tau\in\mc A_*, \sigma\in\mc B_*). \]
Here we identify $(\mc A_* \proten^p \mc B_*)'$ with $\mc{CB}_p(\mc B_*,\mc A)$,
instead of $\mc{CB}_p(\mc A_*,\mc B)$, for convenience, as above we have
been working mainly with the map $R$, and not $L$.  The other choice
follows by symmetry, of course.

\begin{proposition}\label{delta_map_contractive}
With notation as above, and giving $\mc A_*$ and $\mc B_*$ the dual
structures, we have that $\delta$ is a $p$-complete contraction.
\end{proposition}
\begin{proof}
Let $T\in\mathbb M_n(\mc A\vnten\mc B)$, let $\sigma\in\mathbb M_m(\mc B_*)$,
and let $a = ((\delta)_n(T))_m(\sigma) \in \mathbb M_{n\times m}(\mc A)$.
Notice that
\[ a_{ik,jl} = \delta(T_{ij})(\sigma_{kl}) = L(\sigma_{kl})(T_{ij})
\qquad (1\leq i,j\leq n, 1\leq k,l\leq m). \]
We shall, for the proof, give $\mc A_*$ the quotient structure in order
to evaluate the norm on $\mathbb M_{n\times m}(\mc A)$.  Let $\tau\in
\mathbb M_r(\mc A_*)$, and let $\epsilon>0$.  We may find $\hat\tau\in
\mathbb M_r(\mc N(E))$ such that $\hat\tau$ maps to $\tau$, and
$\|\hat\tau\|_r\leq \|\tau\|_r+\epsilon$.  As in the proof of
Lemma~\ref{slice_map_lemma}, we $p$-completely isometrically identify
$\mc B(E\otimes_p F)$ with $\mc{CB}_p(\mc N(E),\mc B(F))$ by the map $\Lambda$.
Then we have that
\begin{align*}
\| \langle\ip{a}{\tau} \|
&= \big\| \langle\ip{ L(\sigma_{kl})(T_{ij}) }{ \tau_{st} }\rangle \big\|
= \big\| \langle\ip{ R(\tau_{st})(T_{ij}) }{ \sigma_{kl} }\rangle \big\|
= \big\| \langle\ip{ \Lambda(T_{ij})(\hat\tau_{st}) }{ \sigma_{kl} }\rangle \big\| \\
&\leq \big\| ((\Lambda)_n(T))_r(\hat\tau) \big\|_{n\times r} \|\sigma\|_m
\leq \big\| ((\Lambda)_n(T)) \big\|_n \|\hat\tau\|_r \|\sigma\|_m \\
&\leq \|\Lambda\|_{pcb} \|T\|_n \|\hat\tau\|_r \|\sigma\|_m
\leq \|T\|_n \|\sigma\|_m (\|\tau\|_r+\epsilon).
\end{align*}
As $\tau$ was arbitrary, we see that $\|a\|_{n\times m} \leq \|T\|_n
\|\sigma\|_m$.  As $\sigma$ was arbitrary, we see that $\|(\delta)_n(T)\|_{pcb}
\leq \|T\|_n$.  Finally, as $T$ was arbitrary, we conclude that $\delta$
is a $p$-complete contraction, as required.
\end{proof}

Now give $\mc A_*$ and $\mc B_*$ the quotient structures.  Then by
Proposition~\ref{proj_is_proj}, the obvious map
\[ \pi_* : \mc N(E)\proten^p \mc N(F) \rightarrow \mc A_*
\proten^p \mc B_* \]
is a $p$-complete quotient map.  Thus
\[ \pi := \pi_*' : (\mc A_*\proten^p\mc B_*)'\rightarrow
\mc B(E\otimes_p F) \]
is a $p$-complete isometry.

\begin{theorem}\label{alg_case_nice}
With notation as above, the map $\pi$ is a weak$^*$-homeomorphic
$p$-completely isometric map with range equal to $\mc A \otimes_{\mc F} \mc B$.
Furthermore, $\pi$ takes $\overline{\mc A\otimes\mc B}$, defined to
be the weak$^*$-closure of $\mc A\otimes\mc B$ in $(\mc A_*\proten^p
\mc B_*)'$, onto $\mc A\vnten\mc B$.
\end{theorem}
\begin{proof}
This follows as for operator spaces (given properties of $\proten^p$
which we established in Proposition~\ref{proj_is_proj}), see
\cite[Proposition~7.2.3]{ER}.
\end{proof}

Finally, we study maps on algebras, and links to complete boundedness.

\begin{theorem}\label{cb_alg_maps}
Let $\phi_1$ and $\phi_2$ be measures, and let $E=L_p(\phi_1)$
and $F=L_p(\phi_2)$.  Let $\mc A\subseteq\mc B(E)$ be a weak$^*$-closed
algebra, and let $M\in\mc{CB}_{p}(\mc A)$ be weak$^*$-continuous.
For any weak$^*$-closed algebra $\mc B\subseteq\mc B(F)$, there exists
a weak$^*$-continuous map $\hat M\in\mc B(\mc A\vnten\mc B)$ such that
$\hat M(a\otimes b) = M(a)\otimes b$ for $a\in\mc A$ and $b\in\mc B$,
and $\|\hat M\| \leq \|M\|_{pcb}$.
\end{theorem}
\begin{proof}
We may suppose that $\phi_1 = \phi_2$ is the counting measure
on $\mathbb N$.  The general case will follow in the same way as
Proposition~\ref{proten_nuc} follows from Proposition~\ref{nuc_p_proj_prod}.
Hence $E=F=\ell_p$.

Let $P_n:\ell_p\rightarrow\ell_p^n$ be the projection onto the first $n$
coordinates, and $\iota_n:\ell_p^n \rightarrow\ell_p$ be the canonical inclusion
map.  Define $\alpha_n:\mc B(\ell_p(\mathbb N\times\mathbb N))
\rightarrow \mc B(\ell_p) \vnten \mc B(\ell_p^n) = \mathbb M_n(\mc B(\ell_p))$ by
\[ \alpha_n(T) = \big( L(P_n'(\delta_i^*)\otimes\iota_n(\delta_j))(T) \big)_{ij}
\qquad (T\in\mc B(\ell_p(\mathbb N\times\mathbb N)), 1\leq i,j\leq n). \]
Let $x=(x_j)_{j=1}^n\subseteq\ell_p$ and $\mu=(\mu_i)_{i=1}^n\subseteq\ell_{p'}$,
and define $y = \sum_{j=1}^n x_j\otimes\iota_n(\delta_j) \in \ell_p\otimes_p\ell_p$
and $\lambda = \sum_{i=1}^n \mu_i\otimes P'_n(\delta_i^*)$.  Then
\[ \|y\| = \Big\| (I\otimes\iota_n)\Big( \sum_{j=1}^n x_j\otimes\delta_j \Big) \Big\|
\leq \Big( \sum_{j=1}^n \|x_j\|^p \Big)^{1/p}, \]
and similarly $\|\lambda\|^{p'} \leq \sum_i \|\mu_i\|^{p'}$.  Then
\begin{align*}
&|\ip{\mu}{\alpha_n(T)(x)}| = \Big| \sum_{i,j=1}^n
   \ip{\mu_i}{L(P'_n(\delta_i^*)\otimes\iota_n(\delta_j))(T)(x_j)} \Big| \\
&= \Big| \sum_{i,j=1}^n \ip{T}{\big( (\mu_i\otimes x_j)
   \otimes P'_n(\delta_i^*)\otimes\iota_n(\delta_j) \big)} \Big|
= \Big| \ip{T}{ \sum_{i,j=1}^n (\mu_i\otimes P'_n(\delta_i^*))
   \otimes (x_j\otimes\iota_n(\delta_j)) } \Big| \\
&= |\ip{T}{\lambda\otimes y}| \leq \|T\|
   \Big( \sum_{i=1}^n \|\mu_i\|^{p'} \Big)^{1/p'}
   \Big( \sum_{j=1}^n \|x_j\|^{p} \Big)^{1/p}.
\end{align*}
Thus $\|\alpha(T)\| \leq \|T\|$, so that $\alpha$ is a contraction.
It is easy to show that $\alpha$ is weak$^*$-continuous.
We have defined $\alpha$ in such a way that $\alpha(T\otimes S)
= T \otimes P_n S\iota_n$ for $S,T\in\mc B(\ell_p)$.

In a similar way, we may define a weak$^*$-continuous contraction
$\beta:\mathbb M_n(\mc B(\ell_p)) \rightarrow \mc B(\ell_p(\mathbb N\times
\mathbb N))$ such that $\beta(T\otimes S) = T\otimes\iota_n SP_n$ for
$T\in\mc B(\ell_p)$ and $S\in\mathbb M_n$.

By weak$^*$-continuity, we see that $\alpha_n(T)\in\mathbb M_n(\mc A)$
for $T\in\mc A\vnten\mc B(\ell_p)$.  As $M\in\mc{CB}_p(\mc A)$, by definition,
we have that $(M\otimes I_n)\alpha_n:\mc A\vnten\mc B(\ell_p)\rightarrow
\mathbb M_n(\mc A)$ is bounded, with $\|(M\otimes I_n)\alpha_n\| \leq \|M\|_{pcb}$.
Thus $\beta_n(M\otimes I_n)\alpha_n:\mc A\vnten\mc B(\ell_p) \rightarrow
\mc A\vnten\mc B(\ell_p)$ is bounded with $\|\beta_n(M\otimes I_n)\alpha_n\|\leq
\|M\|_{pcb}$.  As $\alpha_n, \beta_n$ and $M\otimes I_n$ are weak$^*$-continuous,
so is $\beta_n(M\otimes I_n)\alpha_n$.

Let $(n_\alpha)$ be a subnet of $\mathbb N$ such that the net
$\beta_{n_\alpha}(M\otimes I_{n_\alpha})\alpha_{n_\alpha}(T)$ converges
in the weak$^*$-topology, for each $T\in\mc A\vnten\mc B(\ell_p)$,
say converging to $M_0(T)\in\mc A\vnten\mc B(\ell_p)$.  Then $M_0$ is linear
and bounded, with $\|M_0\| \leq \|M\|_{pcb}$.  Then, for
$i,j,k,l\in\mathbb N$, $a\in\mc A$ and $S\in\mc B(\ell_p)$,
\begin{align*}
& \lim_\alpha \ip{\delta_i^*\otimes\delta_j^*}
   {\beta_{n_\alpha}(M\otimes I_{n_\alpha})\alpha_{n_\alpha}
   (a\otimes S)(\delta_k\otimes\delta_l)} \\
&= \lim_\alpha \ip{\delta_i^*}{M(a)(\delta_k)}
   \ip{\delta_j^*}{\iota_{n_\alpha}P_{n_\alpha}S\iota_{n_\alpha}P_{n_\alpha}(\delta_l)} \\
&= \ip{\delta_i^*}{M(a)(\delta_k)} \ip{\delta_j^*}{S(\delta_l)},
\end{align*}
as eventually, $\iota_{n_\alpha}P_{n_\alpha}(\delta_l)=\delta_l$ and so forth.
Thus 
\[ M_0(a\otimes S) = \lim_\alpha \beta_{n_\alpha}(M\otimes I_{n_\alpha})
\alpha_{n_\alpha}(a\otimes S)  =  M(a)\otimes S
\qquad (a\in\mc A, S\in\mc B(\ell_p)), \]
with the limit taken in the weak$^*$-topology.

Let $\mc A_* = \mc N(\ell_p)/{^\perp}\mc A$ be the predual of $\mc A$,
and let $m\in\mc B(\mc A_*)$ be such that $m'=M$.  Let $\theta:
\mc A_*\otimes\mc N(\ell_p)\rightarrow (\mc A\vnten\mc B(\ell_p))_*
= \mc N(\ell_p(\mathbb N\times\mathbb N)) / {^\perp}(\mc A\otimes\mc B(\ell_p))$
be the canonical map given by
\[ \ip{a\otimes S}{\theta(\tau\otimes\sigma)} = \ip{a}{\tau}\ip{S}{\sigma}
\qquad (a\in\mc A,S\in\mc B(\ell_p),\tau\in\mc A_*,\sigma\in\mc N(\ell_p)). \]
Then $\theta$ is injective, and we claim that $\theta$ has dense range.
If not, then there exists a non-zero $T\in\mc A\vnten\mc B(\ell_p)$ such
that $\ip{T}{\theta(\tau\otimes\sigma)}=0$ for $\tau\in\mc A_*$ and
$\sigma\in\mc N(\ell_p)$.  There hence exists $x\in\ell_p\otimes_p\ell_p$
and $\mu\in\ell_{p'}\otimes_{p'}\ell_{p'}$ with $\ip{\mu}{T(x)}\not=0$.
By approximation, we may suppose that $x=\sum_{n=1}^N x_n\otimes y_n$ and
$\mu=\sum_{m=1}^M \mu_m\otimes\lambda_m$.  Then
\[ 0 \not= \sum_{n,m} \ip{\mu_m\otimes\lambda_m}{T(x_n\otimes y_n)}
= \sum_{n,m} \ip{T}{\theta\big((\mu_m\otimes x_n+{^\perp}\mc A) \otimes
(\lambda_m\otimes y_n)\big)}, \]
a contradiction.
For $a\in\mc A,S\in\mc B(\ell_p),\tau\in\mc A_*$ and $\sigma\in\mc N(\ell_p)$,
we have that
\[ \ip{M_0(a\otimes S)}{\theta(\tau\otimes\sigma)}
= \ip{a}{m(\tau)} \ip{S}{\sigma} = \ip{a\otimes S}{\theta(m(\tau)\otimes\sigma)}. \]
We hence see that $m\otimes I$ extends continuously to a bounded map
on $(\mc A\vnten\mc B(\ell_p))_*$, and so by weak$^*$-density,
$M_0$ is weak$^*$-continuous.

Finally, for $a\in\mc A$ and $b\in\mc B$, we have that
$a\otimes b\in\mc A\vnten\mc B\subseteq\mc A\vnten\mc B(\ell_p)$, and
$M_0(a\otimes b) = M(a)\otimes b$.  As $M_0$ is weak$^*$-continuous,
we hence see that $M_0(\mc A\vnten\mc B) \subseteq \mc A\vnten\mc B$,
and so we may set $\hat M$ to be $M_0$ restricted to $\mc A\vnten\mc B$,
completing the proof.
\end{proof}

\section{Fig\`a-Talamanca-Herz algebras}

We shall briefly introduce the Fig\`a-Talamanca-Herz algebras,
following the notation of \cite{Herz1} (which means that, compared
to some authors, we swap the indexes $p$ and $p'$).

Let $G$ be a locally compact group, and let $\lambda_p:G\rightarrow
\mc B(L_p(G))$ be the \emph{left regular representation}, defined by
\[ \lambda_p(s)(f)(t) = f(s^{-1}t) \qquad (s,t\in G, f\in L_p(G)). \]
We shall also need to use the \emph{right regular representation},
which is defined by
\[ \rho_p(s)(f)(t) = f(ts) \Delta_G(s)^{1/p} \qquad (s,t\in G,f \in L_p(G)), \]
where $\Delta_G$ is the modular function of $G$.
See Section~\ref{multisec} for further details about group
representations.  Let $C(G)$ be the space of continuous functions
from $G$ to $\mathbb C$, let $C_{00}(G)\subseteq C(G)$ be the subspace
of functions with compact support, and let $C_0(G)$ be its closure.
We then define a map $\Lambda_p:L_{p'}(G) \proten L_p(G) \rightarrow C_0(G)$ by
\[ \Lambda_p(g\otimes f)(s) = \ip{g}{\lambda_p(s)(f)}
\qquad (s\in G, f\in L_p(G), g\in L_{p'}(G)). \]
That $\Lambda_p$ maps into $C(G)$ follows as $\lambda_p$ is continuous;
that $\Lambda_p$ maps into $C_0(G)$ follows as $C_{00}(G)$ is
dense in $L_p(G)$ and $L_{p'}(G)$.
Then $A_p(G)$ is defined to be the \emph{coimage} of $\Lambda_p$.
That is, we identify the image of $\Lambda_p$ with the Banach space
$L_{p'}(G)\proten L_p(G) / \ker \Lambda_p$, the latter defining the
norm on $A_p(G)$.  As shown in \cite{Herz1}, $A_p(G)$ becomes
a Banach algebra under pointwise operations.  When $p=2$, $A_2(G)$
agrees with the Fourier Algebra $A(G)$, as studied in \cite{Eymard}.

By standard Banach space results, we see that the dual of $A_p(G)$
may be identified with the space
\[ PM_p(G) = \{ T\in\mc B(L_p(G)) : \ip{T}{\tau}=0
\ (\tau \in \ker\Lambda_p) \}. \]
Notice that $\lambda_p(G) = \{ \lambda_p(s):s\in G\} \subseteq
PM_p(G)$, and that the weak$^*$-closure of $\lambda_p(G)$ is equal
to $PM_p(G)$.  It is then easy to show that $PM_p(G)$ is a
subalgebra of $\mc B(L_p(G))$ (see, for example, \cite[Section~10]{Pier}).
When $p=2$, we have that $PM_2(G) = VN(G)$, the group von Neumann
algebra of $G$.  The duality between $A_p(G)$ and $PM_p(G)$ is
\[ \ip{T}{\Lambda_p(g\otimes f)} = \ip{g}{T(f)}
\qquad (T\in PM_p(G), g\in L_{p'}(G), f\in L_p(G)). \]

As $PM_p(G)\subseteq\mc B(L_p(G))$, we see that $PM_p(G)$ carries
a natural $p$-operator space structure.  As in Section~\ref{weakstaralgs},
we may hence induce the dual $p$-operator space structure on $A_p(G)$.
Alternatively, we may induce the quotient structure on $A_p(G)$,
by defining the map $\Lambda_p: \mc N(L_p(G))\rightarrow A_p(G)$
to be a $p$-complete quotient map.

When $G$ is amenable, the algebra $PM_p(G)$ is easier to handle.
In particular, we have \cite[Theorem~5]{Herz2}, which shows that
when $G$ is amenable, we have that $PM_p(G) = CONV_p(G)
:= \{ T\in\mc B(E) : T\rho_p(s) = \rho_p(s)T \ (s\in G) \}$ 

\begin{theorem}\label{when_amen_thm}
Let $G$ be an amenable locally compact group.  Then the two natural
$p$-operator space structures on $A_p(G)$, as defined above, agree.
\end{theorem}
\begin{proof}
Let $E = L_p(G)$.
By Proposition~\ref{proj_implies_struc_equal}, it suffices to prove
that there is a $p$-completely contractive projection from $\mc B(E)$
onto $PM_p(G)$.  We shall now show that when $G$ is amenable,
such a projection exists.

We may $E$ as a left $L_1(G)$-module for the $\rho_p$ action, so that
$L_1(G)_E^c = CONV_p(G)$, in the notation of Section~\ref{amen_ba}.
Combining Theorem~\ref{amen_char} and Proposition~\ref{quasi_exp}
yields that there is a contractive projection $\mc Q:\mc B(E)
\rightarrow CONV_p(G) = PM_p(G)$.

However, we need to show that $\mc Q$ is actually $p$-completely
contractive.  Let $(d_\alpha)$ be an approximate diagonal of bound one
for $L_1(G)$, and let $d_\alpha = \sum_{n=1}^\infty a^{(\alpha)}_n
\otimes b^{(\alpha)}_n\in L_1(G)\proten L_1(G)$ for each $\alpha$.
Let $T\in\mathbb M_n(\mc B(E))$, let $(x_j)_{j=1}^n\subseteq E$,
and let $(\mu_i)_{i=1}^n\subseteq E'$, so that
\begin{align*} \Big|\sum_{i,j=1}^n & \ip{\mu_i}{\mc Q(T_{ij})(x_j)}\Big|
= \lim_\alpha \Big|\sum_{i,j=1}^n \sum_k \ip{\mu_i}{\rho_p(a^{(\alpha)}_k)
T_{ij}\rho_p(b^{(\alpha)}_k)(x_j)} \Big| \\
&\leq \lim_\alpha \sum_k \Big|\sum_{i,j=1}^n \ip{\rho_p(a^{(\alpha)}_k)'(\mu_i)}
   {T_{ij}\rho_p(b^{(\alpha)}_k)(x_j)} \Big| \\
&\leq \lim_\alpha \sum_k \|T\|_n
   \Big( \sum_{i=1}^n \|\rho_p(a^{(\alpha)}_k)'(\mu_i)\|^{p'} \Big)^{1/p'}
   \Big( \sum_{j=1}^n \|\rho_p(b^{(\alpha)}_k)(x_j)\|^p \Big)^{1/p} \\
&\leq \|T\|_n \lim_\alpha \sum_k \|\rho_p(a^{(\alpha)}_k)\| \|\rho_p(b^{(\alpha)}_k)\|
   \Big( \sum_{i=1}^n \|\mu_i\|^{p'} \Big)^{1/p'}
   \Big( \sum_{j=1}^n \|x_j\|^p \Big)^{1/p} \\
&\leq \|T\|_n \Big( \sum_{i=1}^n \|\mu_i\|^{p'} \Big)^{1/p'}
   \Big( \sum_{j=1}^n \|x_j\|^p \Big)^{1/p}.
\end{align*}
Thus $\|\mc Q\|_n\leq 1$, and so $\|\mc Q\|_{pcb}=1$, as required.
\end{proof}

In \cite[Section~1.31]{Paterson}, the class of groups $G$ such that $PF_2(G)$
is an amenable Banach algebra is discussed: it is somewhat larger
than the class of amenable groups.  When $PF_2(G)$ is amenable,
by weak$^*$-density, we see that $VN(G) = PM_2(G)$ is Connes-amenable,
and this is enough to ensure a projection $\mc B(L_2(G))$ to
$PM_2(G)$ (actually, such a projection is automatically completely
positive, and hence completely contractive, see \cite[Chapter~XV,
Corollary~1.3]{tak3}).
For example, \cite[Page~84]{Paterson} shows that $VN(SL(2,\mathbb R))$
is Connes-amenable, while $SL(2,\mathbb R)$ is not amenable.
Of course, in the $p=2$ case the above theorem is not necessary.
In the $p\not=2$ case, we are not aware of a systematic investigation of when
$PM_p(G)$, for $p\not=2$, is Connes-amenable (see \cite[Theorem~4.4.13]{RundeBook}
for some partial results).  Furthermore, even if we have a projection
$\mc B(L_p(G))\rightarrow PM_p(G)$, it is unclear that this projection
is necessarily $p$-completely contractive.  It seems possible that
the above proof could hence be extended to some non-amenable groups.

However, the existence of a projection onto $PM_p(G)$ is very far from
being necessary, so it also seems possible that another method
of proof could extend the above result to a much larger class of
groups (or even maybe all groups).

We know that $p$-operator spaces are much easier to work with when
they embed into an $L_p$ space.  Henceforth, we shall assume that
$A_p(G)$ carries the dual structure.  We shall resort to the above theorem
when it is necessary to use the quotient structure (which is in many ways
the more natural structure).

Our next task is to show that $A_p(G)$ is an algebra is the category
of $p$-operator spaces.  This is equivalent to saying that the algebra
product defines a bounded (indeed, contractive) map $\Delta:A_p(G) \proten^p
A_p(G) \rightarrow A_p(G)$.  Suppose that $\Delta':PM_p(G)\rightarrow
(A_p(G)\proten A_p(G))' = \mc{CB}_p(A_p(G),A_p(G)') = 
\mc{CB}_p(A_p(G),PM_p(G))$ is a $p$-complete contraction.  Then so
is $\Delta''$, and hence also $\Delta'' \kappa_{A_p(G)\proten^p A_p(G)}
= \kappa_{A_p(G)} \Delta$.  As $\kappa_{A_p(G)}$ is a $p$-complete
isometry, we conclude that $\Delta$ is a $p$-complete contraction.

Define $PM_p(G)\vnten PM_p(G) \subseteq \mc B(L_p(G\times G))$,
as in Section~\ref{algebras}.

\begin{proposition}\label{weak_star_ten_prod}
Let $G$ and $H$ be locally compact groups.  Then $PM_p(G)\vnten
PM_p(H) = PM_p(G\times H)$.
\end{proposition}
\begin{proof}
By definition, $PM_p(G)\vnten PM_p(H)$ is the weak$^*$-closure
of $PM_p(G)\otimes PM_p(G)$ in $\mc B(L_p(G)\otimes_p L_p(H))
= \mc B(L_p(G\times H))$.  For this proof, let $\lambda_p^G:G
\rightarrow\mc B(L_p(G))$ be the left-regular representation,
and define $\lambda_p^H$ and $\lambda_p^{G\times H}$ similarly.
Then it is simple to verify that
\[ \lambda_p^G(s) \otimes \lambda_p^H(t)
= \lambda_p^{G\times H}(s,t) \qquad (s\in G, t\in H). \]
Hence we see immediately that $PM_p(G\times H) \subseteq
PM_p(G) \vnten PM_p(H)$, as $PM_p(G\times H)$ is the weak$^*$-closure
of the span of the image of $\lambda_p^{G\times H}$.

Conversely, we shall show that $\lambda^G_p(G)\otimes PM_p(H)
\subseteq PM_p(G\times H)$, and by symmetry that $PM_p(G)
\otimes \lambda^H_p(H) \subseteq PM_p(G\times H)$.  Thus,
for $S\in PM_p(G)$ and $T\in PM_p(H)$, we have that
\[ S\otimes T = (S\otimes \lambda_p^H(e_H))
( \lambda_p^G(e_G)\otimes T) \in PM_p(G\times H), \]
where $e_G$, $e_H$, is the unit of $G$, respectively $H$.
As $PM_p(G\times H)$ is weak$^*$-closed, we conclude that
$PM_p(G)\vnten PM_p(H) \subseteq PM_p(G\times H)$, completing the proof.

To show that $\lambda^G_p(G)\otimes PM_p(H) \subseteq PM_p(G\times H)$
we shall show that $\ker\Lambda_p^{G\times H} \subseteq
{^\perp}( \lambda_p^G(G) \otimes PM_p(H))$.  Let $\tau =
\sum_{n=1}^\infty \mu_n\otimes x_n \in \ker\Lambda_p^{G\times H}
\subseteq L_{p'}(G\times H) \otimes L_p(G\times H)$.  We regard
$L_p(G\times H)$ as $L_p(G,L_p(H))$, and so we regard each $x_n$
as a function from $G$ to $L_p(H)$.  Similarly $L_{p'}(G\times H)
= L_{p'}(G,L_{p'}(H))$.  Fix $u\in G$, so that
\[ 0 = \ip{\lambda_p^G(u)\otimes\lambda_p^H(v)}{\tau}
= \sum_{n=1}^\infty \int_G \ip{\mu_n(s)}{\lambda_p^H(v)(x_n(u^{-1}s))}
\ ds \qquad ( v\in H ). \]
For each $n$, define $y_n\in L_p(G,L_p(H))$ by $y_n(s) =
x_n(u^{-1}s)$ for $s\in G$.  Thus
\[ 0 = \sum_{n=1}^\infty \ip{\mu_n}{(I\otimes\lambda_p^H(v))(y_n)}
\qquad (v\in H). \]
By using Herz's ideas in \cite[Lemma~0]{Herz1}, this implies that
\[ 0 = \sum_{n=1}^\infty \ip{\mu_n}{(I\otimes T)(y_n)}
= \ip{\lambda_p^G(u)\otimes T}{\tau} \qquad (T\in PM_p(H)). \]
As $u\in G$ was arbitrary, the proof is complete.
\end{proof}

Define $W:L_p(G\times G)\rightarrow L_p(G)$ by
\[ (Wf)(s,t) = f(s,st) \qquad (f\in L_p(G\times G), s,t\in G), \]
so that $W$ is an invertible isometry.  
Define
\[ \Gamma: PM_p(G) \rightarrow PM_p(G)\vnten PM_p(G);
\quad T \mapsto W^{-1}(T\otimes I) W \qquad (T\in PM_p(G)). \]
Let $f\in L_p(G\times G)$ and $s\in G$.  Then
\begin{align*} \big( \Gamma(\lambda_p(s))(f) \big)(r,t) &=
\big( W^{-1} (\lambda_p(s)\otimes I) W(f) \big)(r,t) \\
&= \big( (\lambda_p(s)\otimes I) W(f) \big)(r,r^{-1}t) \\
&= (Wf)(s^{-1}r,r^{-1}t) = f(s^{-1}r,s^{-1}t), \end{align*}
for $r,t\in G$.  Thus $\Gamma(\lambda_p(s)) = \lambda_p(s)\otimes
\lambda_p(s)$.

Recall the definition of the map $\delta:PM_p(G) \vnten PM_p(G)
\rightarrow (A_p(G)\proten^p A_p(G))'$, which is a $p$-complete
contractive by Proposition~\ref{delta_map_contractive}.
For $a,b\in A_p(G)$ and $s\in G$, we have that
\[ \ip{\delta\Gamma(\lambda_p(s))}{a\otimes b} = \ip{\lambda_p(s)\otimes
\lambda_p(s)}{a\otimes b} = a(s) b(s) = (ab)(s)
= \ip{\lambda_p(s)}{\Delta(a\otimes b)}. \]
Thus $\Delta' = \delta\Gamma$.  In particular, as $\Gamma$ is clearly a
$p$-complete contraction, so is $\Delta'$, as required.

\begin{theorem}\label{tenprod_herzalg}
Let $G$ and $H$ be amenable locally compact groups.  Then
$A_p(G) \proten^p A_p(H) = A_p(G\times H)$.
\end{theorem}
\begin{proof}
This proof is an adaptation of \cite[Theorem~7.2.4]{ER}.
By Theorem~\ref{when_amen_thm}, we have that the two $p$-operator
space structures agree on $A_p(G)$ and $A_p(H)$.  By
Theorem~\ref{alg_case_nice}, the map $\pi_*:\mc N(L_p(G))
\proten^p \mc N(L_p(H)) \rightarrow A_p(G)\proten^p A_p(H)$
is a $p$-complete quotient map, so that $\pi=\pi_*':(A_p(G)
\proten^p A_p(H))' \rightarrow \mc B(L_p(G\times H))$ is a $p$-complete
isometry onto its range, which is $PM_p(G) \otimes_{\mc F} PM_p(H)$.

For $w\in\mc N(L_p(G))$, recall the definition of $R(w)$ from
Section~\ref{algebras}.  Let $T\in PM_p(G) \otimes_{\mc F} PM_p(H)
\subseteq \mc B(L_p(G)\otimes_p L_p(H))$, so by definition
$R(w)(T)\in PM_p(H)=CONV_p(H)$ for each $w\in\mc N(L_p(G))$.  Thus,
for $s\in H$, $R(w)(T) \rho_p(s) = \rho_p(s) R(w)(T)$.
By weak$^*$-continuity, this implies that
\[ R(w)\big( T(I\otimes\rho_p(s)) \big)
= R(w)\big( (I\otimes\rho_p(s))T \big). \]
As $w$ is arbitrary, this is that $T(I\otimes\rho_p(s))
= (I\otimes\rho_p(s))(T)$ for each $s\in H$.  By symmetry, we also
see that $(\rho_p(t)\otimes I)T = T(\rho_p(t)\otimes I)$ for $t\in G$.
Consequently $T$ commutes with $\rho_p((t,s))$ for $(t,s)\in G\times H$,
that is, $T\in CONV_p(G\times H) = PM_p(G\times H)$, as $G\times H$
is amenable.

Thus $PM_p(G) \otimes_{\mc F} PM_p(H) = PM_p(G) \vnten PM_p(H) =
PM_p(G\times H)$.  As $\pi$ is a homeomorphism, we conclude
that $(A_p(G) \proten^p A_p(H))' = A_p(G\times H)' = PM_p(G\times H)$
$p$-completely isometrically.  As the quotient and dual structures
agree on $A_p(G\times H)$, and $\pi=\pi_*'$ is weak$^*$-continuous,
this implies that $A_p(G) \proten^p A_p(H) = A_p(G\times H)$, as required.
\end{proof}

In the above proof we use the fact that when $G$ is an amenable group,
we have that $PM_p(G) = CONV_p(G)$.  As communicated to us by
Professor Fig\`a-Talamanca, in \cite{cowling}, Cowling shows that
$PM_p(G) = CONV_p(G)$ for $G=SL(2,\mathbb R)$ and $G=\mathbb F_2$.  Actually,
the proof for $\mathbb F_2$ is not correct, but can be corrected using
results of Haagerup, as done in \cite[Theorem~4.9, Chapter~8]{FTP}.
It is apparently unknown if $PM_p(G) = CONV_p(G)$ for all groups $G$.
We conclude that the main sticking point in this section is
Theorem~\ref{when_amen_thm}.

Finally, we shall show that $A_p(G)$ is amenable in the category of
$p$-operator spaces if and only if $G$ is an amenable group.  By
``amenable in the category of $p$-operator spaces'', we mean that every
$p$-completely bounded derivation from $A_p(G)$ to a $p$-completely
contractive dual $A_p(G)$-bimodule is inner.  The equivalence of this
to $A_p(G)$ having an approximate diagonal in $A_p(G) \proten^p A_p(G)$
follows from exactly the same argument as used for amenability of
Banach algebras (compare with \cite[Section~2]{Ruan1}).  We shall make
heavy use of the already established result in the $p=2$ case,
which is \cite[Theorem~3.6]{Ruan1}.

\begin{theorem}
Let $G$ be a locally compact group.  Then $A_p(G)$ is $p$-operator
space amenable if and only if $G$ is an amenable group.
\end{theorem}
\begin{proof}
Suppose that $A_p(G)$ is $p$-operator space amenable.  Then,
in particular, $A_p(G)$ has a bounded approximate identity, and
so by Leptin's Theorem (compare \cite[Theorem~6]{Herz2}) $G$ is amenable.

Conversely, suppose that $G$ is an amenable group.  Then $A_p(G)
\proten^p A_p(G) = A_p(G\times G)$.  As $G\times G$ is amenable,
by \cite[Theorem~C]{Herz1}, identification of functions gives
a norm-decreasing homomorphism $A_2(G\times G)\rightarrow A_p(G\times G)$
which has dense range.  By Ruan's Theorem, $A_2(G\times G)=A_2(G)
\proten^2 A_2(G)$ contains a bounded approximate diagonal, and hence
so does $A_p(G\times G)$.  Thus $A_p(G)$ is $p$-operator space amenable.
\end{proof}

\subsection{Further homological properties}\label{Fur_homo_prop}

Amenability fits into the study of Hochschild cohomology of Banach
algebras, and there are further (co)homological properties of
Banach algebras which are widely studied.  See
\cite[Chapter~4]{RundeBook} for an introduction to these ideas.
As for amenability, when $A(G)$ is considered as an operator space,
homological properties of $A(G)$ depend upon the group $G$ in the
same (or dual) way to the way that properties of $L_1(G)$ depend upon $G$.

In \cite{Wood}, Wood considers \emph{biprojectivity}, and shows
that $A(G)$ is biprojective (with the operator space structure)
if and only if $G$ is discrete.  Conversely, Helemskii (see
\cite{Hele}) showed that $L_1(G)$ is biprojective (as a Banach
alegbra) if and only if $G$ is compact (and we view discreteness
and compactness as being dual properties, as in the abelian case).

First, some terminology.  Let $\mc A$ be a Banach algebra, let
$E$ and $F$ be $\mc A$-bimodules, and let $\theta\in\mc B(E,F)$.
We say that $\theta$ is an \emph{module homomorphism} if
$\theta(a\cdot x\cdot b) = a\cdot\theta(x)\cdot b$ for $a,b\in\mc A$
and $x\in E$.  We say that $\theta$ is \emph{admissible} if there
exists $\phi\in\mc B(F,E)$ with $\theta\phi\theta = \theta$.
We say that an $\mc A$-bimodule $E$ is \emph{biprojective} when, given
$\mc A$-bimodules $F$ and $G$, a surjective, admissible module
map $\phi:F\rightarrow G$ and a module map $\theta:E\rightarrow G$,
there exists a module map $\psi:E\rightarrow F$ with $\phi\psi = \theta$.

In \cite{Wood}, Wood first adapts these ideas to the category
of operator spaces.  Subject to some technicalities (as usual, to
do with duality) it seems rather likely that this carries over
easily to the $p$-operator space situtation.  Wood next proves that
the multiplication map $A(G) \proten^2 A(G) \rightarrow A(G)$ is
surjective.  Thus uses a number of results, including that
$A(G)\proten^2 A(G) = A(G\times G)$ for all groups $G$, which we
have not been able to generalise to the $A_p(G)$ case.  Furthermore,
this fact is again used in the proof of the main theorem,
\cite[Theorem~4.5]{Wood}.

A Banach algebra $\mc A$ is \emph{weakly-amenable} when every
bounded derivation to $\mc A'$ is inner.  When $\mc A$ is commutative,
this is equiavlent to the (more natural) condition that every
derivation into a \emph{symmetric} $\mc A$-bimodule $E$ is zero.
Here an $\mc A$-bimodule $E$ is \emph{symmetric} if $a\cdot x=
x\cdot a$ for each $a\in\mc A$ and $x\in E$.  It is easy to
translate these conditions into the category of operator spaces,
and in \cite{Spronk} Spronk shows that $A(G)$ is always weakly-amenable
in the category of operator spaces.

Again, we can translate these ideas over to $p$-operator spaces,
but, again, we find that we need properties of the projective
tensor norm which we have not been able to establish in full
generality (it is, of course, pointless to restrict to amenable
groups $G$, as then $A_p(G)$ is amenable, and so trivially weakly-amenable).
Furthermore, Spronk uses simple facts about representations on Hilbert
spaces which seem unlikely to hold for $SQ_p$ spaces, as we lack things
like orthogonal projections.  ??It would be interesting to see if a
different approach could be found to deal with the $p\not=2$ case??

In \cite{samei}, Samei develops the theory of algebras he called
\emph{hyper-Tauberian}, and uses this theory to give a simple and
elegant proof that $A(G)$ is weakly-amenable, as an operator space.
Indeed, Samei's argument easily extends to the $A_p(G)$ algebras,
when given the operator space structure constructed in \cite{LNR}.
This operator space structure suffers from the same issue we have,
in that $A_p(G) \proten^2 A_p(G)$ need not, semmingly, be anything
useful, when $G$ is not amenable.  Samei sidesteps this issue by
first working with $A(G)$ and then transfering the result to $A_p(G)$,
see \cite[Theorem~28]{samei} for details.

We can immediately adapt Samei's definition of what it means to
be hyper-Tauberian to the $p$-operator space setting, and show that
a hyper-Tauberian algebra is weakly-amenable.  It remains to show
that $A_p(G)$ is indeed hyper-Tauberian as a $p$-operator space.
Unfortunately, we again hit a problem here, as we cannot lift
results from $A(G)$ to $A_p(G)$ (as $A(G)$ is not a $p$-operator space!)
and a direct argument, at least following Samei, would again require
us to know what $A_p(G) \proten^p A_p(G)$ is.  It at least seems possible
that a new direct argument could work for $A_p(G)$ in the $p$-operator
space setting, but we have not been able to make progress in this direction.


\section{Multipliers}\label{multisec}

In this section we shall study multipliers of Fig\`a-Talamanca-Herz Algebras.
Much of the hard work is already in the literature, but often without
direct connections being drawn.  We try to collect together these
results in a unified setting here.

It shall be helpful to sketch some results on group representations.
Let $G$ be a locally compact group, and let $E$ be a \emph{reflexive}
Banach space.  We shall define a \emph{group representation of
$G$ on $E$} to be a group homomorphism $\pi:G\rightarrow\mc B(E)$
such that $\pi(s)$ is an isometry for each $s\in G$, and for each
$x\in E$ and $\mu\in E'$, the map $G\rightarrow\mathbb C; s\mapsto
\ip{\mu}{\pi(s)(x)}$ is continuous.  Then $\pi$ extends to a norm-decreasing
homomorphism $\pi:L_1(G)\rightarrow\mc B(E)$ by integration.

We show now sketch the converse to this, which is folklore.  Let
$\pi:L_1(G)\rightarrow\mc B(E)$ be a norm-decreasing homomorphism.
As is standard (see \cite[Theorem~3.3.23]{Dales} for example)
$L_1(G)$ contains an approximate identity $(e_\alpha)$
of bound $1$.  For $s\in G$ and $f\in L_1(G)$, define $s\cdot f
\in L_1(G)$ by $(s\cdot f)(t) = f(s^{-1}t)$ for $t\in G$.  We
may define a map $\sigma:G\rightarrow\mc B(E)$ by
\[ \ip{\mu}{\sigma(s)(x)} = \lim_\alpha \ip{\mu}{(s\cdot e_\alpha)(x)}
\qquad (x\in E,\mu\in E'). \]
Then there exists a subspace $F$ of $E$ such that,
by restriction, $\sigma$ becomes a group representation
$\sigma:G\rightarrow\mc B(F)$.  In fact, there is a contractive
projection $P:E\rightarrow F$ such that $P\pi(f)P = \pi(f)$
for $f\in L_1(G)$, so that the action of $\pi$ on the kernel
of $P$ is trivial, and so we loose nothing by restricting to $F$.
Applying the previous paragraph to $\sigma$ yields the homomorphism
$\pi$, restricted to $F$.  By the Cohen Factorisation Theorem, we have
that $F= \{ \pi(f)(x) : x\in E, f\in L_1(G) \}$.

Now define a map $\Pi:E'\proten E\rightarrow C(G)$ by
\[ \Pi(\mu\otimes x)(s) = \ip{\mu}{\pi(s)(x)}
\qquad (\mu\otimes x\in E'\proten E, s\in G). \]
Here $C(G)$ is the space of continuous functions on $G$; that $\Pi$
maps into $C(G)$ follows by the continuity assumption on $\pi$.
We let $A(\pi)$ be the co-image of $\Pi$: that is, $A(\pi)$ is the
image of $\Pi$ in $C(G)$, but with the norm induced by identifying
$A(\pi)$ with the quotient $E'\proten E / \ker\Pi$.
As explained by Herz in \cite{Herz1}, the obvious definition of
equivalent group representations is a rather strong condition, while
$A(\pi)$ gives a more interesting notion of equivalence (for example,
$A(\pi)$ is one-dimensional if and only if $\pi$ is trivial).

Recall the left-regular representation $\lambda_p:G\rightarrow\mc B(L_p(G))$.
Then $A_p(G) = A(\lambda_p)$.  Let $\pi:G\rightarrow\mc B(E)$ be
some group representation, and let $I_E:G\rightarrow\mc B(E)$ be the
trivial representation on $E$.  Herz shows that $\lambda_p$ has the
useful property that $A(\lambda_p\otimes\pi) = A(\lambda_p\otimes I_E)$
(this is also referred to as \emph{Fell's absorption principle}).
Furthermore, if $E\in SQ_p$ (or, in Herz's terminology, $E$ is a
\emph{$p$-space}) then $A(\lambda_p\otimes I_E) = A(\lambda_p)$
(\cite[Lemma~0]{Herz1}).

For a commutative Banach algebra $\mc A$, we say that a linear map
$T:\mc A\rightarrow\mc A$ is a \emph{multiplier}, denoted by $T\in
\mc M(\mc A)$, if $T(ab) = aT(b)$ for $a,b\in\mc A$.  Then $\mc M(\mc A)$
becomes a Banach algebra with respect to the operator norm.  For a locally
compact group $G$, using the fact that $A_p(G)$ is a \emph{regular
tauberian algebra} (see \cite[Section~3]{Herz2}), we may use the
Closed Graph Theorem to show that each multiplier on $A_p(G)$ is
bounded, and furthermore, each multiplier is given by pointwise
multiplication by some (necessarily continuous) function $u:G\rightarrow
\mathbb C$.  Henceforth we shall treat $\mc M(A_p(G))$ as a subspace
of $C(G)$, with the norm
\[ \|u\|_{\mc M} = \sup\{ \|ua\|_{A_p} : a\in A_p(G), \|a\|_{A_p}\leq 1 \}
\qquad (u\in\mc M(A_p(G))). \]

It is common in the literature to write $B_p(G)$ for $\mc M(A_p(G))$.
This is confusing, as it is standard to denote by $B(G)$ the
\emph{Fourier-Stieltjes Algebra} of $G$.  However, by results of
Nebbia and Losert (see \cite{los}) we have that $B_2(G)=B(G)$
if and only if $G$ is amenable (see \cite[Page~187]{Paterson} for an
example where this confusion arises).  To further confuse the issue, Herz
himself defined a space $B_p(G)$ in \cite{Herz3}, using a notion of Schur
multipliers (which we shall study further below).  Finally, Runde
defined a generalisation of $B(G)$ in \cite{Runde2} which he,
reasonably, denotes by $B_p(G)$.  We shall stick to writing $\mc M(A_p(G))$.

In \cite{DH}, De Canni\`{e}re and Haagerup study
\emph{completely bounded multipliers} of $A_p(G)$, denoted by
$\mc M_0(A_p(G))$.  We have that $B_2(G) = \mc M_0(A_p(G))$ in
Herz's notation (see \cite{BF} where unpublished results of J.\ Gilbert
are used to show this).  Similar ideas are explored \cite{jol}.
We use \cite{DH} and \cite{jol} to motivate to make the following
definitions.

\begin{definition}
Let $G$ be a locally compact group, let $1<p<\infty$, and let
$u\in\mc M(A_p(G))$.  Then $u\in\mc M_{cb}(A_p(G))$ if and only if
$u$ defines a member of $\mc{CB}_p(A_p(G))$ where $A_p(G)$ is given
the dual $p$-operator space structure.  We give $\mc M_{cb}(A_p(G))$
the $p$-completely bounded norm.

We define $M_0(A_p(G))$ to be the space of those functions
$u:G\rightarrow\mathbb C$ such that there exists $E\in SQ_p$ and bounded
continuous maps $\alpha:G\rightarrow E$ and $\beta:G\rightarrow E'$
such that $u(ts^{-1}) = \ip{\beta(t)}{\alpha(s)}$ for $s,t\in G$.
We give $\mc M_0(A_p(G))$ the obvious norm.
\end{definition}

Then, for example, Jolissaint shows in \cite{jol} that $\mc M_0(A_2(G))
= \mc M_{cb}(A_2(G))$.

\begin{lemma}
Let $G$ be a locally compact group, and let $u:G\rightarrow\mathbb C$
be a function.  Then the following are equivalent:
\begin{enumerate}
\item $u\in\mc M(A_p(G))$;
\item There exists a bounded, weak$^*$-continuous operator
   $M:PM_p(G)\rightarrow PM_p(G)$ such that
   $M(\lambda_p(s)) = u(s)\lambda_p(s)$ for $s\in G$.
\end{enumerate}
\end{lemma}
\begin{proof}
Suppose that (1) holds, let $m\in\mc B(A_p(G))$ be the
operator defined by pointwise multiplication by $u$, and let
$M=m'\in\mc B(PM_p(G))$.  Then obviously
$M(\lambda_p(s)) = u(s)\lambda_p(s)$ for $s\in G$.

Conversely, if (2) holds, then as $M$ is weak$^*$-continuous,
there exists $m\in\mc B(A_p(G))$ with $m'=M$.  For $a\in A_p(G)$,
we then have that
\[ u(s)a(s) = \ip{M(\lambda_p(s))}{a} =
\ip{\lambda_p(s)}{m(a)} = m(a)(s) \qquad (s\in G), \]
so that $u$ is pointwise multiplication by $u$, and hence
$u\in\mc M(A_p(G))$.
\end{proof}

When $p=2$ the above can be significantly improved, essentially
because $A(G)$ is a \emph{closed} ideal in $B(G)$; see
\cite[Proposition~1.2]{DH}.

\begin{theorem}
Let $G$ be a locally compact group, and let $1<p<\infty$.  Then $\mc M_0(A_p(G))$
and $\mc M_{cb}(A_p(G))$ are commutative Banach algebras.
Furthermore, $\mc M_{cb}(A_p(G)) = \mc M_0(A_p(G))$ isometrically, and
$\mc M_{cb}(A_p(G)) \subseteq \mc M(A_p(G))$ contractively.
\end{theorem}
\begin{proof}
For the proof, write $\mc M$ for $\mc M(A_p(G))$ and so forth.
Obviously $\mc M_{cb} \subseteq \mc M$ contractively,
from which it follows easily that $\mc M_{cb}$ is a commutative
Banach algebra.  For $E,F\in SQ_p$, by considering the space $E\oplus F$ with
the $\ell_p$ norm $\|(x,y)\| = (\|x\|^p + \|y\|^p)^{1/p}$ for $x\in E,y\in F$,
it follows that $\mc M_0$ is a vector space.  Similarly, by
considering the infinite $\ell_p$ sum of a countable family $(E_n)_{n=1}^\infty
\subseteq SQ_p$, it follows that $\mc M_0$ is a Banach space.  Finally,
by using a suitable tensor product construction for $SQ_p$ spaces
(see \cite[Section~3]{Runde2}) it follows that $\mc M_0$ is a commutative
Banach algebra.

Now let $u\in\mc M_0$ be defined by $u(ts^{-1}) = \ip{\beta(t)}{\alpha(s)}$,
using some $E\in SQ_p$.  Let $x\in L_p(G)$ and $\mu\in L_{p'}(G)$, and let
$a=\Lambda_p(\mu\otimes x)\in A_p(G)$.  Define $\hat x\in L_p(G,E) = L_p(G)
\otimes_p E$ and $\hat\mu\in L_{p'}(G,E')$ by
\[ \hat x(s) = x(s) \alpha(s^{-1}), \quad
\hat\mu(s) = \mu(s) \beta(s^{-1}) \qquad (s\in G). \]
Then $\|\hat x\| \leq \|\alpha\|_\infty \|x\|$, $\|\hat\mu\|\leq
\|\beta\|_\infty \|\mu\|$, and for $s\in G$,
\begin{align*} \ip{\hat\mu}{(\lambda_p(s)\otimes I_E)(\hat x)}
&= \int_G \ip{\hat\mu(t)}{\hat x(s^{-1}t)} \ dt
= \int_G \ip{\beta(t^{-1})}{\alpha(t^{-1}s)} \mu(t) x(s^{-1}t) \ dt \\
&= \int_G u(tt^{-1}s) \mu(t) x(s^{-1}t) \ dt
= u(s) a(s).
\end{align*}
By Herz, we have that $A(\lambda_p\otimes I_E)=A_p(G)$, so that
$ua\in A_p(G)$ and $\|ua\|_{A_p} \leq \|\alpha\|_\infty \|\beta\|_\infty
\|x\| \|\mu\|$.  As $u$ was arbitrary, and by linearity and the
definition of the norms on $A_p(G)$ and $\mc M_0$, we see that
$\mc M_0 \subseteq \mc M$ is a norm-decreasing inclusion.

We can ``amplify'' this argument to show that
$\mc M_0 \subseteq \mc M_{cb}$ contractively.
Let $u\in\mc M_0$ be as before, and let $M\in\mc B(PM_p(G))$ be induced by $u$,
as given by the previous lemma.  Given $x\in L_p(G)$
and $\mu\in L_{p'}(G)$, define $\hat x\in L_p(G,E)$ and $\hat\mu\in
L_{p'}(G,E')$ as above.  It is a simple calculation to show that
\[ \ip{\mu}{M(T)(x)} = \ip{\hat\mu}{(T\otimes I_E)(\hat x)}
\qquad (T\in PM_p(G)). \]
Let $n\in\mathbb N$ and let $T=(T_{ij})\in\mathbb M_n(PM_p(G))$.  Let
$\mu=(\mu_i)_{i=1}^n\in L_{p'}(G)\otimes_p\ell_{p'}^n$ and
$x=(x_j)_{j=1}^n \in L_p(G)\otimes_p\ell_p^n$.  Define $\hat x\in
L_p(G) \otimes_p E \otimes_p \ell_p^n$ by
$\hat x = \sum_{j=1}^n \hat x_j \otimes \delta_j$, so that
\[ \|\hat x\| = \Big( \sum_{j=1}^n \|\hat x_j\|^p \Big)^{1/p}
\leq \|\alpha\|_\infty \Big( \sum_{j=1}^n \|x_j\|^p \Big)^{1/p}
= \|\alpha\|_\infty \|x\|. \]
Similarly define $\hat\mu$, so that $\|\hat\mu\| \leq
\|\beta\|_\infty \|\mu\|$.  Finally, define $S\in
\mc B(L_p(G)\otimes_p E\otimes_p\ell_p^n)$ by
\[ S(x\otimes y\otimes\delta_j) = \sum_{i=1}^n
T_{ij}(x) \otimes y \otimes \delta_i
\qquad (x\in L_p(G), y\in E, 1\leq j\leq n). \]
If $\phi: L_p(G)\otimes_p E\otimes_p\ell_p^n \rightarrow
L_p(G)\otimes_p\ell_p^n\otimes_p E$ is the canonical isometry,
then $\phi S \phi^{-1} = T\otimes I_E$, so that $\|S\|=\|T\|$.
Then
\begin{align*}
|\ip{\mu}{(M)_n(T)(x)}| &=
\Big| \sum_{i,j=1}^n \ip{\mu_i}{M(T_{ij})(x_j)} \Big|
= \Big| \sum_{i,j=1}^n \ip{\hat\mu_i}{(T_{ij}\otimes I_E)(\hat x_j)} \Big| \\
&= | \ip{\hat\mu}{S(\hat x)} | \leq \|\hat\mu\| \|\hat x\| \|S\|
\leq \|\alpha\|_\infty \|\beta\|_\infty \|T\| \|\mu\| \|x\|,
\end{align*}
so that $M\in\mc{CB}_p(PM_p(G))$ with $\|M\|_{pcb}\leq
\|\alpha\|_\infty \|\beta\|_\infty$, as required.

To show that $\mc M_{cb} \subseteq \mc M_0$, one can easily
adapt Jolissaint's proof in \cite{jol} by combining it with Pisier's
representation theorem for $p$-completely bounded maps
(Theorem~\ref{pisier_pcb}), a task we now sketch.  Let $u\in\mc M_{cb}
\subseteq\mc M$, and let $M\in\mc B(PM_p(G))$ be given as in the
lemma above.  By definition, $M\in\mc{CB}_p(PM_p(G))$, so as $PM_p(G)$
is a unital algebra, by the comment after Theorem~\ref{pisier_pcb},
there exists $E\in SQ_p$, a $p$-representation
$\hat\pi:PM_p(G)\rightarrow\mc B(E)$ and $U:L_p(G)\rightarrow E$
and $V:E\rightarrow L_{p'}(G)$ with $\|U\| \|V\|\leq\|M\|_{pcb}$,
such that $M(T) = V\hat\pi(T)U$ for $T\in PM_p(G)$.  It is clear from
the definitions that $\hat\pi$ is a norm-decreasing algebra homomorphism,
and so $\hat\pi\circ\lambda_p:L_1(G)\rightarrow\mc B(E)$ is a
norm-decreasing algebra homomorphism.  By the discussion at the beginning
of this section, there hence exists a one-complemented subspace $F$ of $E$
and a group representation $\sigma:G\rightarrow\mc B(F)$.  As the action
of $\hat\pi\circ\lambda_p$ is only non-trivial of $F$, and $F$ is
one-complemented, we loose nothing by assuming that actually $E=F$.
We then notice that
\[ V \sigma(s) U = u(s) \lambda_p(s) \qquad (s\in G). \]

Choose $\mu_0\in L_{p'}(G)$ and $x_0\in L_p(G)$ with $\|x_0\|=\|\mu_0\|=
\ip{\mu_0}{x_0}=1$, and define $\alpha:G\rightarrow E$ and
$\beta:G\rightarrow E'$ by
\[ \alpha(s) = \sigma(s^{-1})U\lambda_p(s)(x_0), \quad
\beta(s) = \sigma(s)'V'\lambda_p(s^{-1})'(\mu_0) \qquad (s\in G), \]
so that $\|\alpha\|_\infty \leq \|U\|$ and $\|\beta\|_\infty\leq
\|V\|$.
Hence, for $s,t\in G$, we have that
\begin{align*}
\ip{\beta(t)}{\alpha(s)} &= 
\ip{\sigma(t)'V'\lambda_p(t^{-1})'(\mu_0)}{\sigma(s^{-1})U\lambda_p(s)(x_0)} \\
&= \ip{\mu_0}{\lambda_p(t^{-1})V\sigma(ts^{-1})U\lambda_p(s)(x_0)} \\
&= u(ts^{-1}) \ip{\mu_0}{\lambda_p(t^{-1})\lambda_p(ts^{-1})\lambda_p(s)(x_0)}
= u(ts^{-1}).
\end{align*}
It remains to show that $\alpha$ and $\beta$ are continuous.  However,
this follows immediately, as a weakly-continuous group representation is
strongly continuous.  Thus $\mc M_{cb} \subseteq \mc M_0$ contractively,
completing the proof.
\end{proof}

\subsection{Herz's Multiplier algebras}

We shall now show how these ideas relate to Herz's algebras $B_p(G)$.  To avoid
confusion, we shall write instead $HS_p(G)$, for \emph{Herz-Schur} multiplier.
Let $I$ be an index set, and let $\psi:I\times I\rightarrow\mathbb C$ be
a function.  We say that $\psi\in V_p(I)$ if and only if, for each
$T\in\mc B(\ell_p(I))$, we have that $T\psi\in\mc B(\ell_p(I))$, where
$T\psi$ is defined by
\[ \ip{\delta^*_i}{(T\psi)(\delta_j)} =
\psi(i,j) \ip{\delta^*_i}{T(\delta_j)} \qquad (i,j\in I). \]
By the closed-graph theorem, $V_p(I) \subseteq \mc B(\mc B(\ell_p(I)))$,
which gives the obvious norm on $V_p(I)$.

Let $X$ be a separable locally compact space, and let $X_d$ be the
space $X$ equipped with the discrete topology.  Then we set $V_p(X)$ to
be $C(X\times X) \cap V_p(X_d)$.  Finally, suppose that $G$ is
a separable locally compact group, and let $u\in FS_p(G)$ if and only
if $\psi\in V_p(G)$ where $\psi$ is defined by $\psi(s,t) = u(st^{-1})$
for $s,t\in G$.  For an arbitrary $G$, recall that there is an open
and closed separable subgroup $H$ such that $G$ is the union of left
cosets of $H$.  As such, we can reduce topological questions about $G$ to
questions about $H$, as $G/H$ has the discrete topology.  To avoid tedious
calculations, we shall not mention such topological issues further.

\begin{proposition}\label{Pis_mult_thm}
Let $I$ be an index set, let $\psi:I\times I\rightarrow\mathbb C$
be a function, and let $C>0$.  Then the following are equivalent:
\begin{enumerate}
\item $\psi\in V_p(I)$ and $\|\psi\|_{V_p}\leq C$;
\item There is a measure space $(\Omega,\nu)$ and elements
   $(x_j)_{j\in I}\subseteq L_p(\nu)$ and
   $(\mu_i)_{i\in I}\subseteq L_{p'}(\nu)$ such that
   $\psi(i,j) = \ip{\mu_i}{x_j}$ for each $i,j\in I$, and
   $\sup_i \|\mu_i\| \sup_j\|x_j\|\leq C$;
\item $\psi$ is a $p$-completely bounded multiplier on $\mc B(\ell_p(I))$,
   with $\|\psi\|_{pcb} \leq C$.
\end{enumerate}
\end{proposition}
\begin{proof}
These follow from Theorems~5.11 and~8.2 in \cite{Pisier2}.
\end{proof}

Notice that if $G$ is a discrete group, then using conditions (2)
and (3) above, it is easy to show that $FS_p(G)
= \mc M_0(A_p(G))$ with equal norms.
However, for general $G$, we have the problem that the above
proposition works with $G_d$, hence losing continuity conditions.

Herz shows in \cite[Lemme~1]{Herz3} and \cite[Lemme~2]{Herz3}
that we have the following alternative definition of $V_p(X)$.

\begin{proposition}\label{Alt_Defn_Vp}
Let $X$ be a separable locally compact space, and let $\mu$ be a
Radon measure on $X$ such that each non-empty open subset of $X$
has non-zero $\mu$-measure.  Then $\psi\in V_p(X)$ if and only if
$\psi$ is continuous and there exists $C>0$ such that for
$\mu\in L_{p'}(X,\mu)$ and $x\in L_p(X,\mu)$, there exists
$(\mu_n)_{n=1}^\infty \subseteq L_{p'}(X,\mu)$ and
$(x_n)_{n=1}^\infty \subseteq L_{p}(X,\mu)$ with
\[ \mu(s) x(t) \psi(s,t) = \sum_{n=1}^\infty \mu_n(s) x_n(t)
\qquad (s,t\in X), \]
almost everywhere in $\mu$, with $\sum_{n=1}^\infty \|\mu_n\|_{p'}
\|x_n\|_p \leq C \|\mu\|_{p'} \|x\|_p$.
\end{proposition}

That is, $V_p(X)$ coincides with the space of continuous multipliers
of $L_{p'}(X,\mu) \proten L_p(X,\mu)$, once we have made sense of
what this means.
Let $G$ be a locally compact group with the Haar measure.  Then
the above applies to $V_p(G)$, and hence also to $FS_p(G)$.

Let $G$ be a locally compact group, let $\psi\in V_p(G)$, and let
$n\in\mathbb N$.  Let $G_n = G\times\{1,2,\ldots,n\}$ where
$\{1,2,\ldots,n\}$ is given the counting measure, so that
$L_p(G\times\{1,2,\ldots,n\}) = L_p(G) \otimes_p \ell_p^n$.
Define $\psi_n:G_n\times G_n \rightarrow\mathbb C$ by
\[ \psi_n\big( (s,i) , (t,j) \big) = \psi(s,t)
\qquad (s,t\in G, 1\leq i,j\leq n), \]
so that $\psi_n$ is continuous.  We shall now show that $\psi_n
\in V_p(G_n)$, using the original definition of $V_p$.  Let
$T\in\mc B(\ell_p(G_n))$, so we may also view $T$ as a member
of $\mathbb M_n(\mc B(\ell_p(G)))$, say $T=(T_{ij})$, where
\[ \ip{\delta_s^*}{T_{ij}(\delta_t)}
= \ip{\delta_s^*\otimes\delta_i^*}{T(\delta_t\otimes\delta_j)}
\qquad (s,t\in G, 1\leq i,j\leq n). \]
Let $S = \psi_n\cdot T$, so viewing $S\in\mathbb M_n(\mc B(\ell_p(G)))$,
\begin{align*}
\ip{\delta_s^*}{S_{ij}(\delta_t)}
&= \ip{\delta_s^*\otimes\delta_i^*}{(\psi_n\cdot T)(\delta_t\otimes\delta_j)} \\
&= \psi_n( (s,i),(t,j) ) \ip{\delta_s^*\otimes\delta_i^*}{T(\delta_t\otimes\delta_j)}
= \psi(s,t) \ip{\delta_s^*}{T_{ij}(\delta_t)},
\end{align*}
for $s,t\in G$ and $1\leq i,j\leq n$.  By Proposition~\ref{Pis_mult_thm},
as $\psi$ is automatically $p$-completely bounded, we see that
$\psi_n\in V_p(G_n)$ with $\|\psi_n\|_{V_p} \leq \|\psi\|_{V_p}$.

Now let $u\in FS_p(G)$, so that when $\psi(s,t)=u(st^{-1})$ for
$s,t\in G$, we have that $\psi\in V_p(G)$.  Let $M_u\in\mc B(
L_{p'}(G)\proten L_p(G))$ be the multiplier defined by $\psi$, using
Herz's alternative definition of $V_p(G)$ as shown in
Proposition~\ref{Alt_Defn_Vp}.  Let $x\in L_p(G)$ and
$\mu\in L_{p'}(G)$, so that $a=\Lambda_p(\mu\otimes x)\in A_p(G)$.  Then
\[ \Lambda_p\big( M_u(\mu\otimes x) \big)(s)
= \int_G u(tt^{-1}s) \mu(t) x(s^{-1}t) \ dt
= u(s) a(s) \qquad (s\in G), \]
so that $M_u$ drops under $\Lambda_p$ to pointwise multiplication of
$A_p(G)$ by $u$.  We hence immediately see that $FS_p(G) \subseteq
\mc M(A_p(G))$ contractively.  Combining this observation with the
previous paragraph, we immediately have the following.

\begin{theorem}
Let $G$ be a locally compact group.  Then $FS_p(G) =
\mc M_{cb}(A_p(G))$ isometrically.
\end{theorem}

\subsection{Algebraic definitions}\label{alg_defs}

In \cite{DH}, a more group-theoretic characterisation of
$\mc M_{cb}(A_p(G))$ is shown, and this is used in \cite{BF} to
show that $FS_2(G) = \mc M_{cb}(A(G))$ (which we generalised above,
using another method).

Given sets $I$ and $J$ and functions $u:I\rightarrow\mathbb C,
v:J\rightarrow\mathbb C$, let $u\times v:I\times J\rightarrow\mathbb C$
be defined by $(u\times v)(i,j) = u(i)v(j)$ for $i\in I$ and $j\in J$.

\begin{proposition}
Let $G$ be a locally compact group, let $1<p<\infty$, and let
$u\in\mc M_{cb}(A_p(G))$.  Then, for every locally compact group $H$,
$u\times 1_H \in\mc M(A_p(G\times H))$ and $\|u\times 1_H\|_{\mc M}
\leq \|u\|_{pcb}$.
\end{proposition}
\begin{proof}
By Proposition~\ref{weak_star_ten_prod},
we know that $PM_p(G)\vnten PM_p(H) = PM_p(G\times H)$.
By the above lemma, there exists a weak$^*$-continuous map
$M\in\mc B(PM_p(G))$ such that $M(\lambda_p(s)) = u(s)\lambda_p(s)$
for $s\in G$.  Again, by the lemma, we wish to show that there exists
a weak$^*$-continuous map $\hat M\in PM_p(G\times H)$, such that
\[ \hat M(T\otimes S) = M(T)\otimes S \qquad
(T\in PM_p(G), S\in PM_p(H)). \]
However, this follows immediately from Theorem~\ref{cb_alg_maps},
which also shows that $\|u\times 1_H\|_{\mc M} = \|\hat M\|
\leq \|M\|_{pcb} = \|u\|_{pcb}$.
\end{proof}

In \cite{DH}, the converse to the above is shown in the case $p=2$.
Furthermore, to check that $u$ is completely-bounded, it suffices to
check that $u\times 1_K\in\mc M(A_p(G\times K))$ in the special case
that $K=SU(2)$.  However, we do not have a simple description of what
$PM_p(SU(2))$ is, unless $p=2$.

\subsection{Multipliers and amenability}

In \cite{Runde2}, Runde suggests a definition of a $p$-generalisation
of the Fourier-Stieltjes algebra, which he denotes by $B_p(G)$.
In what follows, we shall follow the conventions of Herz, which means
that we sometimes swap $p$ with $p'$ as compared to Runde.  We define
$B_p(G) \subseteq C(G)$ to be functions of the form
\[ a(s) = \ip{\mu}{\pi(s)(x)} \qquad (s\in G), \]
where $\pi:G\rightarrow\mc B(E)$ is a representation on some $E\in SQ_p$,
and $x\in E$, $\mu\in E'$.  We set $\|a\|_{B_p} = \inf \{ \|\mu\| \|x\|
\}$ where the infimum runs over all representations.  Runde shows that
$B_p(G)$ is a commutative Banach algebra.  It is immediate that
$B_p(G) \subseteq \mc M_0(A_p(G))$ contractively.

It is shown in \cite[Corollary~5.3]{Runde2} that when $G$ is an
amenable locally compact group, we have that $\mc M(A_p(G)) =
B_p(G)$ isometrically, where $B_p(G)$ is Runde's generalisation of
the Fourier-Stieltjes algebra.  We thus immediately have the following.

\begin{proposition}
Let $G$ be an amenable locally compact group, and let $1<p<\infty$.
Then $B_p(G) = \mc M_{cb}(A_p(G)) = \mc M(A_p(G))$ isometrically.
\end{proposition}

As stated above, Nebbia and Losert (see \cite{los}) show that $\mc M(A(G))
= B(G)$ if and only if $G$ is amenable.  In \cite{Bo}, Bo\.{z}ejko
showed that for a discrete group $G$, $\mc M_{cb}(A(G)) = B(G)$ if
and only if $G$ is amenable.  A key point in the proof is that, as a
Banach space, $B(G)$ has cotype 2.  We conjecture that Runde's algebra
$B_p(G)$ has cotype $\max(p,p')$, but we seem to be rather far from
having the tools available to prove this.

In unpublished lecture notes, \cite{loslec}, Losert shows in full
generality that $\mc M_{cb}(A(G)) = B(G)$ only when $G$ is amenable.
The arguments are very close to those used in \cite{los}, but it
appears that it is not possible to simply take the result of \cite{los}
and directly deduce the corresponding result for $\mc M_{cb}(A(G))$.
Furthermore, Losert's arguments in \cite{los} seem to depend upon the
Hilbert space basis of $A(G)$ much more than Nebbia's and Bo\.{z}ejko's
arguments.  We hence seem to be rather far from being able to show that
$\mc M_{cb}(A_p(G)) = B_p(G)$ only when $G$ is amenable, when $p\not=2$.

\section{Conclusions}

Compared to the operator space structure on $A_p(G)$ considered in
\cite{LNR}, we get a contractive quantised Banach algebra, and not
just a bounded algebra product.  It could also be argued that our
approach is more natural, as $A_p(G)$ is an $L_p$-space generalisation
of $A(G)$, so arguably $L_p$ spaces should be used to define a
quantised structure on $A_p(G)$.  However, our approach seems to
require amenability to be introduced to get the theory to work
perfectly.  We are not aware of anyone considering multipliers in
the framework of \cite{LNR}.  It would be interesting to see
if Herz's ideas appear naturally in that setting, as they do in our setting.

It would be interesting to investigate if Theorem~\ref{tenprod_herzalg}
holds for any non-amenable groups, when $p\not=2$.  Furthermore,
it would be interesting to try to extend the tentative results in
Section~\ref{alg_defs}.  Surely a first step in this programme would be
to study the algebras $PM_p(G)$ for, say, $G=SU(2)$.  Finally, surely
the ideas in Section~\ref{Fur_homo_prop} have scope for further study.

We have hinted that perhaps the definition of a $p$-operator space
is not correct.  To be precise, for operators spaces, we consider not
just a space $E$, but also the spaces $\ell_2^n\otimes E$.  This is
reasonable, as $\ell_2^n$ is (up to isometric isomorphism) the only
$n$-dimensional Hilbert space.  For $p$-operator spaces, we replace
$\ell_2^n$ with $\ell_p^n$, but we have less justification for this,
as there are many $n$-dimensional $SQ_p$ spaces.  Of course, Pisier's
and Le Merdy's results suggest that maybe this is enough, as we do
get an intrinsic characterisation of $SQ_p$ spaces, for example.
A more technical problem here is seemingly we do not have a well-defined
way to define a tensor product of two $SQ_p$ spaces.  In
\cite[Section~3]{Runde2}, Runde shows that given $E,F\in SQ_p$, we
may define a completion of $E\otimes F$ in such a way as to get another
$SQ_p$ space, and with a suitable mapping property holding.  However,
it seems that Runde's construction depends upon the chosen representation
of $E$ and $F$ as subspaces of quotients of $L_p$ spaces.

\section{Acknowledgements}

I would like to thank Christian Le Merdy for useful correspondence,
and for supplying me with the counter-example detailed after
Proposition~\ref{kappa_iso}.  I would like to thank Viktor
Losert for sending me a copy of his lecture notes \cite{loslec}.
I would like to thank Alessandro Fig\`a-Talamanca for giving me the
reference to the rarely cited paper \cite{Fig}.  Finally, I would like
to thank Volker Runde for sending me the references which lead to
Section~\ref{Fur_homo_prop}.

\noindent\emph{Author's Address:}
\parbox[t]{3in}{St. John's College,\\
Oxford,\\
OX1 3JP,\\
United Kingdom.}

\bigskip\noindent\emph{Email:} \texttt{matt.daws@cantab.net}

\end{document}